\documentclass{article} 

\usepackage{amssymb,amsmath,epsfig,graphics,psfrag,graphicx,color,mathpazo}
\usepackage{todonotes}
\usepackage{lineno}
\usepackage{fullpage}

\definecolor{lightblue}{rgb}{0.22,0.45,0.70}
\numberwithin{equation}{section}
\numberwithin{figure}{section}
\newcommand\cero{\boldsymbol{0}}
\newcommand\bB{\mathbf{B}}
\newcommand\bC{\mathbf{C}}

\newcommand\bF{\mathbf{F}}

\newcommand\bI{\mathbf{I}}
\newcommand\bM{\mathbf{M}}
\newcommand\bN{\mathbf{N}}
\newcommand\bP{\mathbf{P}}

\newcommand\bX{\mathbf{X}}

\newcommand\beps{\boldsymbol{\varepsilon}}
\newcommand\ff{\boldsymbol{f}}

\newcommand\bk{\boldsymbol{k}}

\newcommand\nn{\boldsymbol{n}}

\newcommand\bsigma{\boldsymbol{\sigma}}

\newcommand\bu{\boldsymbol{u}}
\newcommand\bv{\boldsymbol{v}}
\newcommand\bw{\boldsymbol{w}}
\newcommand\bx{\boldsymbol{x}}

\newcommand\RR{\mathbb{R}}

\newcommand\qon{\quad\hbox{on }\quad}

\newcommand\bdiv{\mathop{\mathbf{div}}\nolimits}
\newcommand\vdiv{\mathop{\mathrm{div}}\nolimits}
\newcommand\bDiv{\mathop{\mathbf{Div}}\nolimits}
\newcommand\vDiv{\mathop{\mathrm{Div}}\nolimits}
\newcommand\bGrad{\mathop{\mathbf{Grad}}\nolimits}
\newcommand\vGrad{\mathop{\mathrm{Grad}}\nolimits}
\newcommand\tr{\mathop{\mathrm{tr}}\nolimits}
\newcommand\bt{\boldsymbol{t}}

\newcommand\bnabla{\boldsymbol{\nabla}}

\newcommand{\ii}{\mathrm{i}}
\newcommand{\e}{\mathrm{e}}
\newcommand{\ev}{\widehat{\mathrm{\boldsymbol{e}}}}
\newcommand{\konp}{\kappa_1 \omega_1 h_1(m_0) + \kappa_2 h_2(m_0)}
\newcommand{\kon}{k_{\text{on}}}
\newcommand{\dkonp}{\kappa_1 \omega_1 h_1^{\prime}(m_0) + \kappa_2 h_2^{\prime}(m_0)}
\newcommand{\koffp}{( 1 - h_1(m_0) )(\kappa_3 + \kappa_4 b_0 )}
\newcommand{\koff}{k_{\text{off}}}
\newcommand{\dmkoffp}{h_1^{\prime}(m_0)(\kappa_3 + \kappa_4 b_0 )}
\newcommand{\dbkoffp}{( 1 - h_1(m_0) )\kappa_4}
\newcommand{\Ap}{\mu + \lambda}
\newcommand{\Bp}{\rho\phi^2 + \mu k^2}
\newcommand{\Cp}{C_0\phi + \frac{\kappa}{\eta}k^2}
\newcommand{\Amp}{\left( \dkonp \right)( 1 - e_0 ) + \dmkoffp e_0}
\newcommand{\Abp}{\dbkoffp e_0}
\newcommand{\Aep}{\konp + \koffp}
\newcommand{\Ep}{\alpha m_0 \e^{-m_0}}
\newcommand{\Hp}{h_3(m_0) + h_3^{\prime}(m_0)m_0}
\newcommand{\taup}{\left( \lambda + \frac{2\mu}{3} \right) \tau \frac{1 - \zeta m_0^2}{\left( 1 + \zeta m_0^2 \right)^2}}

\newtheorem{mydef}{Definition}[section]
\newtheorem{prop}{Proposition}[section]

\newcommand{\dx}{\,\mbox{d}x}

\allowdisplaybreaks

\title{Coupling chemotaxis and growth poromechanics \\for the modelling of feather primordia patterning}
\author{Nicol\'as A. Barnafi\thanks{Department of Mathematics, Università di Pavia, Via Ferrata 1, 27100 Pavia, Italy. E-mail: \texttt{nicolas.barnafi@unipv.it}}, \and Luis Miguel De Oliveira Vilaca\thanks{Laboratory of Artificial \& Natural Evolution (LANE), Department of Genetics and Evolution, University of Geneva, 4 Boulevard d’Yvoy, 1205 Geneva, Switzerland; and SIB Swiss Institute of Bioinformatics, Geneva,
Switzerland. E-mail: \texttt{LuisMiguel.DeOliveiraVilaca@unige.ch, Michel.Milinkovitch@unige.ch}}, \and Michel C. Milinkovitch\footnotemark[2] \and Ricardo Ruiz-Baier\thanks{School of Mathematics, Monash University, 9 Rainforest Walk, Melbourne 3800 VIC, Australia; and Universidad Adventista de Chile, Casilla 7-D, Chillán, Chile. E-mail: \texttt{ricardo.ruizbaier@monash.edu}}}
\date{\today}

\begin{document}
\maketitle
\begin{abstract}
We propose a new mathematical model for the interaction of skin cell populations with fibroblast growth factor and bone morphogenetic protein, occurring within deformable porous media. The equations for feather primordia pattering are based on the work by K.J. Painter \textit{et al.} [J. Theoret. Biol., \textbf{437} (2018) 225--238]. We perform a linear stability analysis to identify relevant parameters in the coupling mechanisms, focusing in the regime of infinitesimal strains. We also extend the model to the case of nonlinear poroelasticity and include solid growth by means of Lee decompositions of the deformation gradient. We present a few illustrative computational examples in 2D and 3D, and briefly discuss the design of tailored efficient solvers. 
\end{abstract}

\section{Introduction}
Chemotaxis models \cite{keller70} describe the directed movement of cells in response to chemicals (attractants or repellents), and can predict the formation of clustered structures. This process has been observed in a variety of embryogenesis processes, such as gastrulation \cite{yang2002} and feather development \cite{lin09}. In the latter, we refer to the buds that then give origin to feathers as \emph{primordia}, and their origin has been modelled through Keller--Segel type chemotaxis models considering fibroblast growth factor (FGF) and bone morphogenic protein (BMP) \cite{mou11, painter12, painter18}.

Poroelasticity is a mixture model in which a solid phase coexists with (at least) one fluid phase \cite{coussy04}. Biological organs and tissues are naturally porous at the tissue level, as they are composed, for example, of both muscle and blood. This phase separation persists even up to the cellular level, as there is the cytoskeleton and the cytoplasm. For this reason, poroelastic models have become very pervasive in the modelling of soft living tissue, such as in oedema formation \cite{barnafi22,lourenco22}, cardiac perfusion \cite{vuong15,barnafi212}, lung characterisation \cite{berger16}, and brain injury \cite{vilaca20}.

In the context of biologically-oriented problems, experiments have shown that the rheology of cytoplasm within living cells exhibits a poroelastic behaviour \cite{moee13}, and in turn, the composition of cells and the extracellular matrix constitutes an overall poromechanical system. The presence of chemical solutes locally modifies morphoelastic properties  and these processes can be homogenised to obtain macroscopic models of poroelasticity coupled with advection-reaction-diffusion equations (see e.g. \cite{collis17,penta14}). {With this biological basis, the scope of our work is twofold: on one hand, we extend the existing pattern formation models for primordia by considering their interaction with the intracellular space in the outset of growth. On the other hand, we perform a thorough stability analysis to understand the coupling mechanisms in the model and to investigate the conditions that  give rise to pattern formation.}

We have structured the remainder of this paper in the following manner. Section \ref{sec:model} describes the coupled model for poro-mechano-chemical interactions, restricting the presentation to the regime of linear poroelasticity. We give an adimensional form of the governing equations, making precise boundary and initial conditions. In Section ~\ref{sec:stability} we perform a linear stability analysis addressing pattern formation according to Turing instabilities. We separate the discussion in some relevant cases and derive and portray patterning spaces. Section \ref{sec:growth} is devoted to extending the model to the case of nonlinear (finite-strain) poroelasticity and material growth, stating also the coupling with chemotaxis in the undeformed configuration. Some numerical examples are given in Section \ref{sec:results} and we close with a summary and discussion of model extensions in Section \ref{sec:concl}.

\section{A coupled model of linear poroelasticity and chemotaxis}\label{sec:model}
Let us consider a piece of soft material as a porous medium in $\mathbb{R}^d$, $d=2,3$, composed by a mixture of
incompressible grains and interstitial fluid, whose description can be placed in the context of
the classical Biot consolidation problem (see e.g.~\cite{showalter00}). In the absence of gravitational forces, of
body loads, and of
mass sources or sinks, we seek for each time $t\in (0,t_{\mathrm{final}}]$,
 the displacement of the  porous skeleton, $\bu(t):\Omega\to \RR^d$, and the pore pressure of the
fluid, $p(t):\Omega\to\RR$, such that
\begin{linenomath*}\begin{subequations}\begin{align}
\partial_t \bigl(C_0 p+\alpha_{BW}\vdiv\bu\bigr) -\frac{1}{\eta}
\vdiv\{\kappa \nabla p\}  &= 0 & \text{in $\Omega\times(0,t_{\mathrm{final}}]$}, \label{eq:mass}\\
\bsigma & = \bsigma_{\text{poroelast}} + \bsigma_{\text{act}} & \text{in $\Omega\times(0,t_{\mathrm{final}}]$}, \label{eq:sigma}\\
\rho\partial_{tt}\bu-\bdiv\bsigma & = \cero & \text{in $\Omega\times(0,t_{\mathrm{final}}]$}, \label{eq:momentum}
\end{align}\end{subequations}\end{linenomath*}
where 
 $\kappa(\bx)$ is the hydraulic conductivity of the porous medium,  $\rho$ is the
density of the solid material, $\eta$ is the constant viscosity of the
interstitial fluid, $C_0$ is the constrained specific storage coefficient, $\alpha_{BW}$ is the Biot-Willis
consolidation parameter, In \eqref{eq:sigma} we are supposing that the
poromechanical deformations are also actively influenced by microscopic tension generation. A very simple
description is given in terms of active stresses: we assume that the total Cauchy stress
contains a passive and an active component, where 
\begin{linenomath*}\begin{equation}\label{eq:total-stress}
\bsigma_{\text{poroelast}} =  \lambda (\vdiv \bu)\bI +2\mu \beps(\bu) - p\bI,
\end{equation}\end{linenomath*}
and $\bsigma_{\text{act}}$
is specified in \eqref{eq:active}, below. The tensor $\beps(\bu) = \frac{1}{2} (\bnabla \bu + \bnabla \bu^\intercal)$ is that  of infinitesimal strains, $\bI$ denotes the second-order identity tensor, and $\mu,\lambda$ are the Lam\'e constants (shear and dilation moduli)
of the solid structure. Equations \eqref{eq:mass}-\eqref{eq:momentum} represent the conservation
of mass, the constitutive relation, and the conservation of linear momentum, respectively.

In addition, let us consider a modified Patlak--Keller--Segel model for the distribution of chemotactic cell populations
of mesenchymal cells, $m$, epithelium activation state, $e$, fibroblast growth factor (FGF), $f$, and bone
morphogenetic protein (BMP), $b$. The base-line model has been developed in \cite{painter18} and
(after being properly modified to account for the motion of the underlying deformable porous media) it can be summarised as
follows
\begin{linenomath*}\begin{subequations} \begin{align}
	\label{eq:m}
  \partial_t m + \partial_t \bu \cdot \nabla m -
	\vdiv\bigl( D_m \nabla m - \alpha m \exp(-\gamma m) \nabla f \bigr)
&= 0 & \text{in } \Omega\times(0,t_{\mathrm{final}}],\\
	\label{eq:e}
\partial_t e -[ \kappa_1 w(\bx,t) h_1(m)+\kappa_2h_2(m)] (1-e)+[1-h_1(m)](\kappa_3+\kappa_4b)e & =0 & \text{in } \Omega\times(0,t_{\mathrm{final}}],\\
\label{eq:f}
  \partial_t f + \partial_t \bu \cdot \nabla f -
	\vdiv( D_f\, \nabla f )
- \kappa_{F}e + \delta_{F}f + \xi_f \vdiv\bu & = 0 & \text{in } \Omega\times(0,t_{\mathrm{final}}],\\
\label{eq:b}
  \partial_t b + \partial_t \bu \cdot \nabla b -
	\vdiv( D_b\, \nabla b )
- \kappa_{B}h_3(m)m + \delta_{B}b& = 0 & \text{in } \Omega\times(0,t_{\mathrm{final}}],
\end{align}\end{subequations}\end{linenomath*}
where $D_m,D_f,D_b$ are positive definite diffusion matrices,  and the spatio-temporal and
nonlinear coefficients  (in this case, the priming wave responsible for the generation of
spatial patterns and the degree of clustering of mesenchymal cells, respectively) are defined as
\begin{linenomath*}\begin{equation}\label{eq:w-h}
	w(\bx,t)  = \frac{\omega_1}{2}\{1+\tanh(\omega_2[t-x_2/\omega_3])\}, \qquad h_i(m) = m^{P_i}[K_i^{P_i}+m^{P_i}]^{-1},
\end{equation}\end{linenomath*}
where $\kappa_i,\omega_i,P_i,\gamma,\delta_i,\xi_f$ are positive model constants.
Note that the mechano-chemical feedback (the process where
mechanical forces modify the reaction-diffusion effects) is here
assumed only through an additional reaction term in the FGF equation \eqref{eq:f},
depending linearly on volume change. Since $\xi_f>0$, this contribution acts as a local sink of FGF during dilation.

On the other hand, we assume that
the active stress component acts isotropically on the medium (see e.g. \cite{jones12}),
and it depends nonlinearly on the concentration of mesenchymal cells, as proposed for instance in \cite{murray88}
\begin{linenomath*} \begin{equation}\label{eq:active}
\bsigma_{\text{act}} =  \biggl(\lambda+\frac{2\mu}{3}\biggr) \frac{\tau m}{1 + \zeta m^2} \bI ,
\end{equation}\end{linenomath*}
with $\tau$ a model constant to be specified later on, which can be positive (implying that
the function modifies the motion of the medium as a local dilation) or negative (isotropically
distributed compression).

In order to reduce the number of model parameters, we focus on a dimensionless counterpart of
systems \eqref{eq:mass}-\eqref{eq:momentum} and \eqref{eq:m}-\eqref{eq:b}, which can be derived using
the following transformation, suggested in \cite{painter18}
\begin{linenomath*}\begin{gather*}
m = \frac{m^*}{\gamma}, \quad e = e^*, \quad f = \frac{\kappa_F f^*}{\delta_B}, \quad b = \frac{\kappa_B b^*}{\gamma}, \\ 
\bu = \sqrt{\frac{D_b}{\delta_B}}\bu^*, \quad p = p^*, \quad
t = \frac{t^*}{\delta_B}, \quad \bx = \sqrt{\frac{D_b}{\delta_B}}\bx^*,
\end{gather*}\end{linenomath*}
where hereafter the stars are dropped for notational convenience. The
adimensional coupled system is then equipped with appropriate initial data at rest
\begin{linenomath*}
\begin{align}\label{eq:initial}
m(0) = m_{0}, \quad e(0)=e_0=\frac{\kappa_1 w(\bx,0) h_1(m_0)+\kappa_2h_2(m_0)}{\kappa_1 w(\bx,0) h_1(m_0)+\kappa_2h_2(m_0)+
[1-h_1(m_0)](\kappa_3+\kappa_4b)}, \nonumber \\
f(0) = \frac{e_0}{\delta_F}, \quad b(0) = h_3(m_0) m_0,\quad
\bu(0)= \cero, \quad \partial_t\bu(0)= \cero,  \quad p(0) = 0, 
\end{align}\end{linenomath*}
defined in $\Omega$; 
and boundary conditions in the following manner
\begin{linenomath*}\begin{subequations}\begin{align}
\label{bc:m}
\{ D_m\, \nabla m - \alpha m \exp(-\gamma m) \nabla f \} \cdot \nn = D_f\, \nabla f \cdot\nn =  \nabla b  \cdot\nn &= 0 &\text{on $\partial\Omega\times(0,t_{\text{final}}]$},\\
\label{bc:Gamma}
\bu = \bu_\Gamma\quad \text{and} \quad \frac{\kappa}{\eta} \nabla p \cdot\nn &= 0 &\text{on $\Gamma\times(0,t_{\text{final}}]$},\\
\label{bc:Sigma}
\bsigma\nn = \bt \quad\text{and}\quad p&=p_0  &\text{on $\Sigma\times(0,t_{\text{final}}]$},
\end{align}\end{subequations}\end{linenomath*}
where $\nn$ denotes the outward unit normal on the boundary, and  $\partial\Omega = \Gamma\cup\Sigma$  is disjointly split into $\Gamma$ and $\Sigma$
where we prescribe clamped boundaries and zero fluid normal fluxes; and a given traction $\bt$
together with constant fluid pressure $p_0$, respectively.

\section{Linear stability analysis and dispersion relation}\label{sec:stability}
\subsection{Preliminaries.}
We proceed to derive a linear stability analysis for the coupled problem \eqref{eq:mass}-\eqref{eq:active}. As usual the analysis is performed on an infinite domain in $\mathbb{R}^d$, with $d = 2, 3$. 
The first step consists in linearising the poro-mechano-chemical system around a steady state, defined in \eqref{eq:initial}. The linearised dimensionless equations are given by
\begin{linenomath*}\begin{equation}\label{eq:LinSysEDP}\begin{split}
	\partial_t \bigl(C_0 p+\alpha_{BW}\vdiv\bu\bigr) -\frac{1}{\eta}
\vdiv\{\kappa \nabla p\}  &= 0,\\
\rho\partial_{tt}\bu-\bdiv\bigl(\bsigma_{\text{poroelast}}+\bsigma_{\text{act}}^{\text{lin}}\bigr) & = \cero, \\
  \partial_t m - \vdiv\{ D_m \nabla m - \Ep \nabla f \}
&= 0,  \\
\partial_t e - A_m(m_0,e_0) m + A_e(m_0,b_0) e + A_b(m_0,e_0) b &= 0,\\
  \partial_t f - \vdiv\{ D_f\, \nabla f \} - e + \delta_{F}f + \xi_f \vdiv\bu & = 0, \\
  \partial_t b - \Delta b - H_3(m_0) m + \delta_{B} b& = 0,\end{split}\end{equation}\end{linenomath*}
where
\begin{linenomath*}
\[\bsigma_{\text{act}}^{\text{lin}} = \left( \taup \right) m \bI,\]
\end{linenomath*}
and $A_m$, $A_e$, $A_b$, $H_3$ are functions of the steady state values, which will be made precise in Definition \ref{def:LinSysMatNoDim}, below. We look then for solutions of the form $\bu, p, m, e, f, b \propto e^{\ii\bk\cdot\bx + \phi t}$, where $\bk$ is the wave vector (a measure of spatial structure) and $\phi$ is the linear growth factor. By substituting this ansatz on the system \eqref{eq:LinSysEDP}, we get a system of linear equations for the vector $\bw = \left(\bu, p, m, e, f, b\right)^\intercal$, where the associated complex eigenvalues $\phi$, give information on occurrence of instability of the steady state (\textit{i.e.}, pattern formation), when its real component is positive.\footnote{We stress that the present analysis will only
address Turing patterning (Hopf bifurcations or other forms of instability that relate to analysing the imaginary part of $\phi$, are not considered).} In order to specify such system, we collect in Proposition \ref{prop:ImpDerivative} some useful preliminary relations. The proof is
postponed to the Appendix.
\begin{prop}\label{prop:ImpDerivative}
	Let $g,\bv$ be sufficiently regular scalar and vector functions defined by $g = \exp(\ii \bk\cdot\bx + \phi t)$ and $\bv = \bv_0\exp(\ii \bk\cdot\bx + \phi t)$ (with $\bv_0$, a constant vector), respectively. Let us also set $\theta = \vdiv\bv$. Then
\begin{linenomath*}
\begin{subequations}
\begin{gather}
		\nabla f = \ii\bk f, \ \Delta f = \ii^2f\bk\cdot\bk = - k^2 f,\ \partial_t f = \phi f, \ 
		\bnabla \bv = \ii \bv \otimes \bk, \ \vdiv\, \bv = \ii \bv \cdot \bk, \ \partial_t \bv = \phi \bv,
		\label{eq:DerivativeFV} \\
\bdiv \beps(\bv)  = -\frac{(\bv\cdot\bk)\bk + k^2\bv}{2},\  \partial_t \beps(\bv) = \phi \beps(\bv), \ 
		\bdiv (\theta \bI) = - (\bv\cdot\bk)\bk, \ \partial_t \theta = \phi \theta.
		\label{eq:DerivativeEpsilonTheta}
\end{gather}
\end{subequations}
\end{linenomath*}
\end{prop}
\noindent The sought system is defined next, starting from \eqref{eq:LinSysEDP}.

\begin{mydef}\label{def:LinSysMatNoDim}
Let $\bw = \left(\bu, p, m, e, f, b\right)^\intercal\in \mathbb{R}^{d+5}$, $d=2,3$ be the vector of independent variables. Then the associated linear system is given by $\bM\bw= \cero_{d+5}$, where the  matrix $\bM_{(d+5) \times (d+5)}$ adopts the form
\begin{linenomath*}	\begin{equation*}
		\bM = \begin{bmatrix}
			\bM_{11}  & \bM_{12} \\
			\bM_{21} & \bM_{22} \\
		\end{bmatrix},
	\end{equation*}\end{linenomath*}
with the blocks defined as
\begin{linenomath*}	\begin{align*}
		\bM_{11} &= \begin{bmatrix}
			A_1k_1 + B & A_1k_2 & \cdots & A_1k_d & \ii k_1 \\
			A_2k_1 & A_2k_2 + B & \cdots & A_2k_d & \ii k_2 \\
			\vdots & \vdots & \ddots & \vdots & \vdots \\
			A_dk_1 & A_1k_2 & \cdots & A_dk_d + B & \ii k_d \\
			\ii \alpha_{BW} \phi k_1 & \ii \alpha_{BW} \phi k_2 & \cdots & \ii \alpha_{BW} \phi k_d & C \\
		\end{bmatrix},\\
		\bM_{12} &= \begin{bmatrix}
			-\ii \taup k_1 & 0 & 0 & 0 \\
			\vdots & \vdots & \vdots & \vdots \\
			-\ii \taup k_d & 0 & 0 & 0 \\
		\end{bmatrix},\\
		 \bM_{21} &= \begin{bmatrix}
			0 & \cdots & 0 \\
			0 & \cdots & 0 \\
			\ii \xi_f k_1 & \cdots & \ii \xi_f k_d \\
			0 & \cdots & 0 \\
		\end{bmatrix}, \\
		\bM_{22}  &= \begin{bmatrix}
			\phi + D_m k^2 & 0 & -\Ep k^2 & 0 \\
			-A_m & \phi + A_e & 0 & A_b \\
			0 & -1 & \phi + D_f k^2 + \delta_F & 0 \\
			-H_3 & 0 & 0 & \phi + k^2 + 1
		\end{bmatrix},
	\end{align*}\end{linenomath*}
and where the relevant entries are defined as $B=\Bp$, $C=\Cp$, and
\begin{linenomath*}	\begin{align*}
	A_i & = \left(\Ap\right) k_i, \quad \text{for all $i$ in $\{1,\ldots,d\}$, }\\
		A_m &= \Amp, \\
		A_e &= \Aep, \\
		A_b &= \Abp, \quad
		H_3 = \Hp.
	\end{align*}\end{linenomath*}
\end{mydef}

We then proceed to obtain a dispersion relation associated with the characteristic polynomial  $P(\phi)=\det(\bM)$ of the matrix described in Definition \ref{def:LinSysMatNoDim}. We obtain
\begin{linenomath*}\begin{equation}
	P(\phi;k^2) = B(\phi;k^2)^{d-1}\left( P_1(\phi;k^2)+P_2(\phi;k^2)P_3(\phi;k^2) \right),
	\label{eq:phiPolyCompactVersion}
\end{equation}\end{linenomath*}
where $B(\phi;k^2)=\Bp$ is a polynomial with pure imaginary roots. Consequently, it does not have an influence on the stability of the steady state of system \eqref{eq:mass}-\eqref{eq:active}. The polynomials $P_i$, $i=1,2,3$, are given in what follows.
\begin{mydef}\label{def:Poly}
The polynomials conforming the characteristic equation \eqref{eq:phiPolyCompactVersion} are
\begin{linenomath*}\begin{align}
	P_1(\phi;k^2) & = b_3 \phi^3 + b_2 \phi^2 + b_1\phi + b_0,\nonumber\\
	P_2(\phi;k^2) &= a_4 \phi^4 + a_3 \phi^3 + a_2 \phi^2 + a_1\phi + a_0,\label{eq:phiPoly}\\
	P_3(\phi;k^2) & = c_3 \phi^3 + c_2 \phi^2 + c_1\phi + c_0,\nonumber
\end{align}
\end{linenomath*}
where all coefficients in \eqref{eq:phiPoly} adopt the following forms
\begin{linenomath*}\begin{align*}
a_0 &= D_mD_f(\kon+\koff)k^6 + \left( (D_m\delta_B+D_mD_f)(\kon+\koff) - \Ep A_m \right)k^4 \\
	&\qquad + \left( D_m\delta_B(\kon+\koff) + \Ep H_3 A_b - \Ep A_m \right)k^2,\\
	a_1 &= D_mD_fk^6 + \left((D_mD_f+D_m+D_f)(\kon+\koff)+D_m\delta_B+D_mDf\right)k^4 \\
	&\qquad + \left( (D_m\delta_B+\delta_B+D_m+D_f)(\kon+\koff) + D_m\delta_B - \Ep A_m \right)k^2 + \delta_B(\kon+\koff), \\
	a_2 &= \left( D_mD_f+D_m+D_f \right)k^4 + \left( D_m\delta_B + (D_m+D_f+1)(\kon+\koff) + \delta + D_m + D_f \right)k^2 \\
	&\qquad + (\delta_B + 1 )( \kon + \koff ) + \delta_B, \\
a_3 &= \left( D_m + D_f + 1 \right)k^2 + \delta_B + \kon + \koff + 1, \qquad 	a_4 = 1, \\
b_i & = -\left(\xi_f \taup \Ep k^2\right) \widehat{b}_i, \qquad
\widehat{b}_0 = \frac{\kappa}{\eta}\left( \kon + \koff \right)\left(k^2+1\right)k^2, \\
\widehat{b}_1 &= \frac{\kappa}{\eta}k^4 + \left( \left( C_0 + \frac{\kappa}{\eta} \right)\left( \kon + \koff \right) + \frac{\kappa}{\eta} \right)k^2 + C_0\left( \kon + \koff \right), \\
\widehat{b}_2 & = \left( C_0 + \frac{\kappa}{\eta} \right)k^2 + C_0\left( \kon + \koff + 1 \right), \qquad
	\widehat{b}_3 = C_0,\\
c_0 &= \frac{\kappa}{\eta}(2\mu+\lambda)k^4, \qquad c_1 = \left( C_0(2\mu + \lambda)+\alpha_{BW} \right) k^2, \qquad
 c_2 = \rho\frac{\kappa}{\eta}k^2 , \qquad 	c_3 = \rho C_0, \\
	  \kon & =\konp, \qquad \koff=\koffp.
\end{align*}\end{linenomath*}
\end{mydef}

From Definition \ref{def:Poly} we immediately see that the characteristic polynomial is of up to seventh degree, making it difficult to determine analytically the main features of the system. 

Next we concentrate on addressing some specific scenarios of particular interest. Unless specified otherwise, we will employ the parameter values provided in Table \ref{tab:LAparams}.
Note that the term $\bsigma_{\text{act}}^{\text{lin}}$ vanishes for $m_0=1$, so, differently from \cite{painter18},
we modify the steady state for mesenchymal cells concentration to $m_0=2$.

\begin{table}[t]
\begin{center}
{\small\begin{tabular}{|l|}
\hline
\hline
$m_0 = 2$, $D_m = 0.01$, $D_f = 0.1$, $\alpha = 4$, $\kappa_1 = 0.05$, $\kappa_2 = 0.025$, $\kappa_3 = 1$, $\kappa_4 = 1$, \\[1.5ex]
$K_1 = 1$, $K_2 = 2$, $K_3 = 5$, $P_1 = P_2 = P_3 = 2$, $\delta_F = 1$, $\omega_1 = 1$,  $\omega_2 = 5$, $\omega_3 = 0.04$\\[1.5ex]
$E = 3\cdot10^4$, $\nu = 0.4$, $C_0 = 10^{-3}$, $\kappa = 10^{-4}$, $\alpha_{BW} = 0.1$, $\tau = 60$, $\eta = 0.1$, $\xi_f = 0.15$, $\rho = 1$, $\zeta = 1$ \\
\hline
\hline
\end{tabular}}
\caption{Model parameters used in the linear stability analysis of the poro-mechano-chemical system.}\label{tab:LAparams}
\end{center}
\end{table}

\subsection{Spatially homogeneous distributions.} If $k^2=0$ then the characteristic polynomial $P(\phi;k^2)$ is
\begin{linenomath*}
\begin{align*}
	P(\phi;0) = & \rho C_0 \phi^8 \bigl[ \phi^3 + \left( \delta_B + \kon + \koff + 1 \right)\phi^2 \\
	&\qquad \quad + \left((\delta_B+1)(\kon+\koff) + \delta_B \right)\phi + \delta(\kon+\koff) \bigr].
\end{align*}
\end{linenomath*}
The well-known Routh--Hurwitz conditions (see e.g. \cite{routh1877}) state that for any polynomial of order 3, a necessary and sufficient set of conditions should be satisfied in order to ensure that all the roots are in the space of complex non-positive real values, $\mathbb{C}_{-}=\{z\in \mathbb{C} \, : \, \Re({z}) \leq 0 \}$.  For a general polynomial $P(\phi)=\alpha_3\phi^3+\alpha_2\phi^2+\alpha_1\phi+\alpha_0$, where all the $\alpha_i>0$, we have
\begin{linenomath*}
\begin{equation}
	\alpha_2\alpha_1 - \alpha_3\alpha_0 > 0.
	\label{eq:RHPoly3}
\end{equation}
\end{linenomath*}
In our case, $\alpha_2\alpha_1 - \alpha_3\alpha_0 = \left( ( \delta_B + 1 )( \kon + \koff ) + \delta_B \right)\left( \delta_B + \kon + \koff \right) + \delta_B + \kon + \koff$, which is a real positive constant. This indicates that the steady state is stable in a spatially homogeneous
system, irrespective of the value of the model parameters.

\subsection{Zero chemotaxis.} When $\alpha=0$, all the coefficients $b_i$ from Definition \ref{def:Poly} are zero. Consequently, the characteristic polynomial is defined by the sub-polynomials $P_2$ and $P_3$ only. As all the parameters implied in $c_i$ are positive and the Routh--Hurwitz condition \eqref{eq:RHPoly3} is satisfied, we know that the $P_3$ polynomial has only non-positive real part complex solutions and therefore does not influence the stability of the steady state. Thus, the instability condition can only come from $P_2$, related only to the chemotaxis model analysed in \cite{painter18}. Substituting $\alpha=0$ in the sub-matrix $\bM_{22}$ leads to a dispersion relation
\begin{linenomath*}
\[P(\phi)=(\phi + D_m k^2)(\phi + A_e)(\phi + D_f k^2 + \delta_F)(\phi + k^2 + 1),\] 
\end{linenomath*}
and consequently the roots are negative real numbers. The absence of chemotaxis prevents then any pattern formation irrespective of the presence or absence of poro-mechanical coupling.

\subsection{Uncoupled system.} Similarly as in the no-chemotaxis scenario above, imposing $\tau = 0$ or $\xi_f=0$ leads to $P_1(\phi;k^2)=0$ for all $\phi, k^2$. Additionally, as $P_3(\phi;k^2)$ has only non-positive roots, only the chemotaxis sub-polynomial enables us to find eigenvalues that lead to unstable systems. For $P_2$ the Routh--Hurwitz conditions are given by
\begin{linenomath*}
\begin{equation}
	a_i > 0, \quad a_2 a_3 - a_1 a_4 > 0 \, \text{ and } \, a_1 a_2 a_3 - a_0 a_3^2 - a_1^2 a_4 > 0.
	\label{eq:RHPoly4}
\end{equation}
\end{linenomath*}
Then, from Definition \ref{def:Poly}, we can readily observe that $a_4$, $a_3$ and $a_2$ are real positive coefficients for any positive parameter  value. A complete analysis of \eqref{eq:RHPoly4} is analytically quite involved. Nevertheless, and in accordance with
the analysis in \cite{painter18}, we can restrict the discussion to the conditions that break $a_0 > 0$.

We first observe, from its definition, that $a_0$ is a polynomial of even order with respect to $k^2$, guaranteeing that it has at least one local extrema. Consequently, we look for the critical wave number, $k^2_c>0$, defined by the equation $a_0^{\prime}(k^2)=0$, that substituting in $a_0$ will lead exactly to be at zero for the associated critical parameter $\theta_c$ that we want to analyse (\textit{i.e.}, $a_0(\theta_c;k_c^2)=0$). In the uncoupled scenario, $a_0^{\prime}$ is a quadratic polynomial with respect to $k^2$ and roots can be easily obtained. Nonetheless, in order to have $a_0$ negative for positive $k^2$, one has to satisfy the following three conditions
\begin{linenomath*}
\begin{subequations}
\begin{align}
	( D_m \delta_B + D_m D_f )( \kon + \koff ) - \Ep A_m & < 0, \label{eq:condUC1}\\
	 D_m\delta_B(\kon+\koff) + \Ep H_3 A_b - \Ep A_m  & < 0, \label{eq:condUC2}\\
	 4\left( (D_m\delta_B+D_mD_f)(\kon+\koff) - \Ep A_m \right)^2 - \qquad \qquad \qquad \qquad \qquad \qquad &  \nonumber \\
	12 D_mD_f(\kon+\koff)\left( D_m\delta_B(\kon+\koff) + \Ep H_3 A_b - \Ep A_m \right) & > 0.
	\label{eq:condUC3}
\end{align}
\end{subequations}
\end{linenomath*}

\begin{figure}[!t]
\begin{center}
\includegraphics[width=0.46\textwidth]{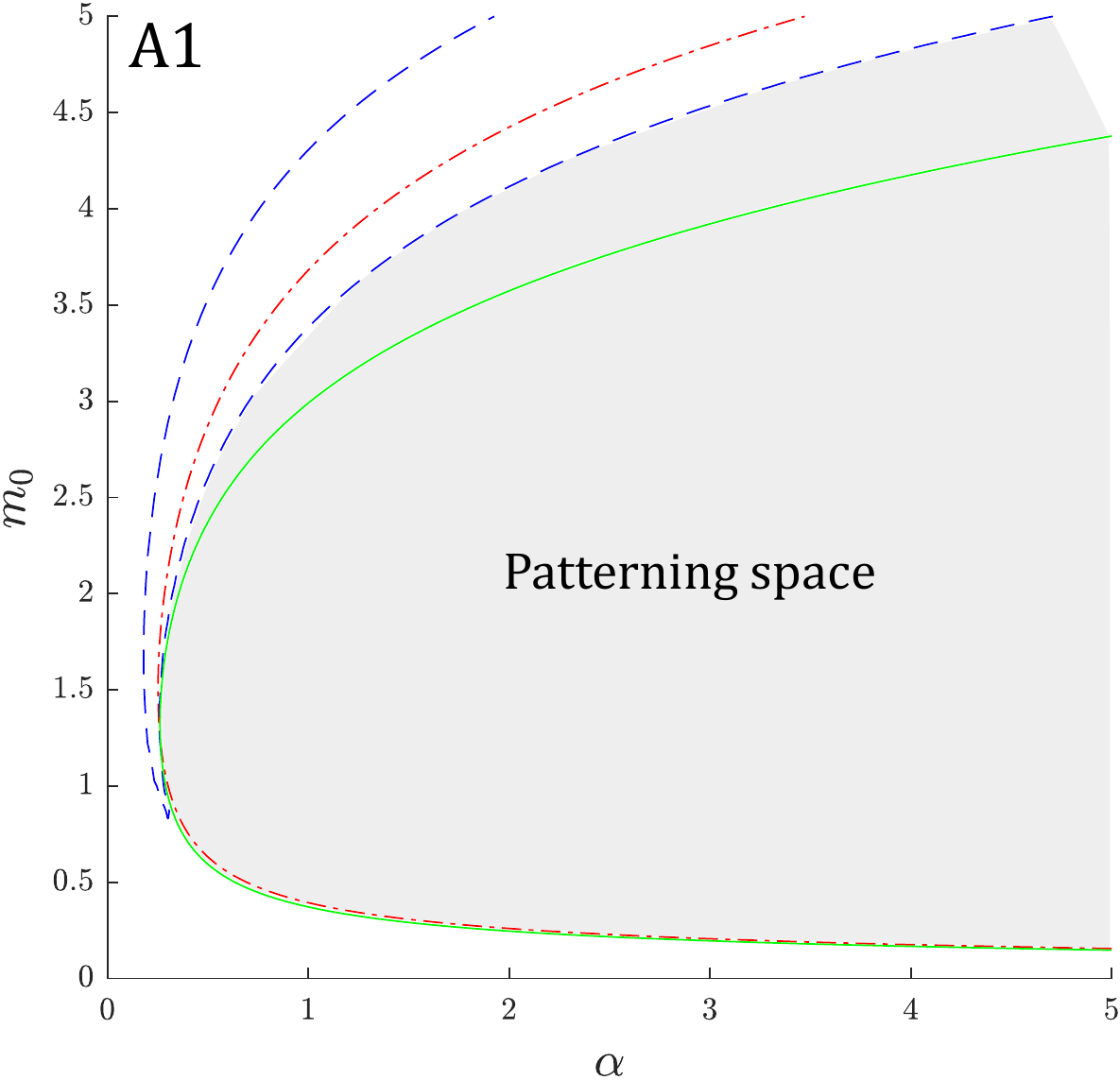}
\includegraphics[width=0.46\textwidth]{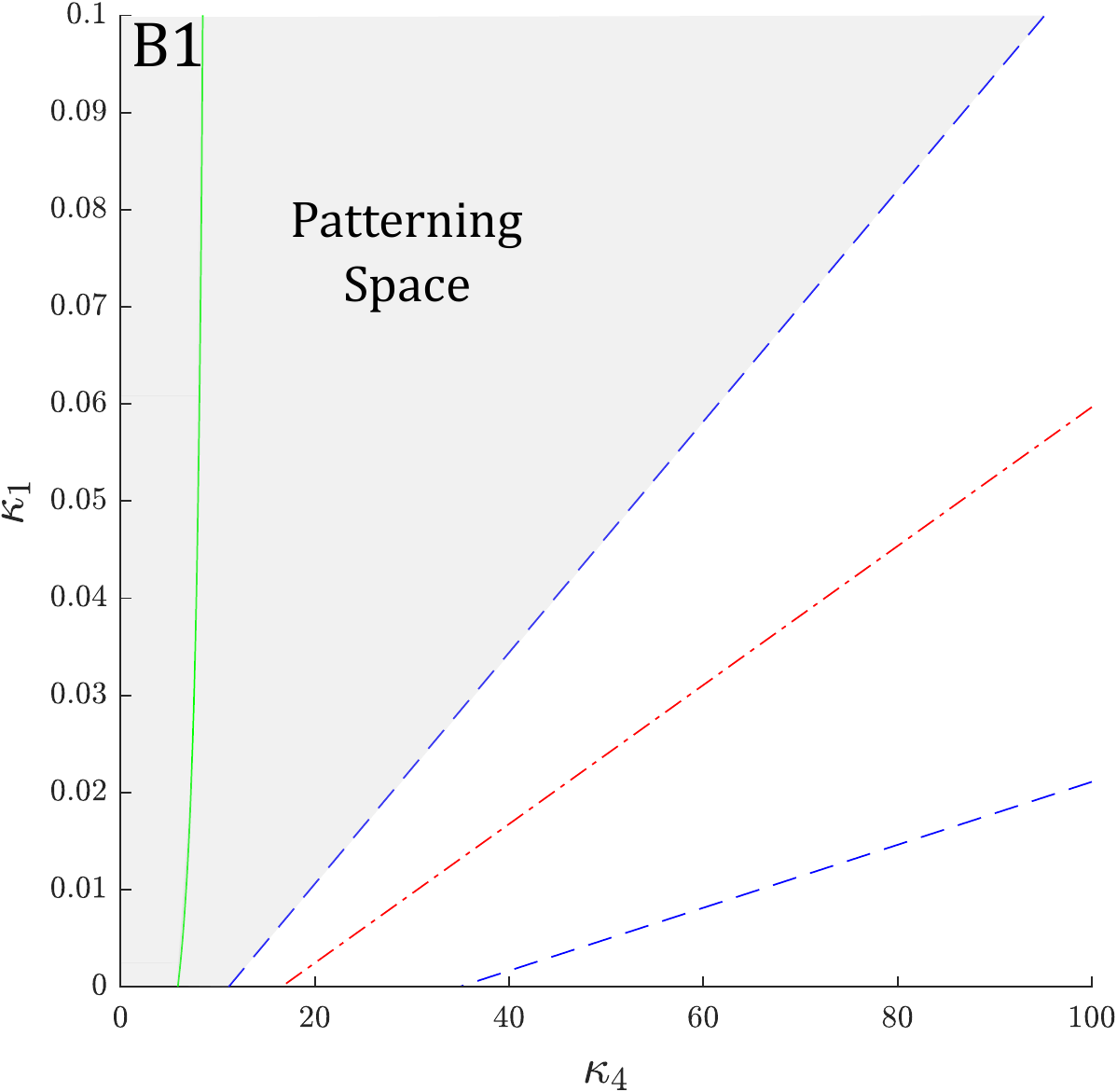}\\
\includegraphics[width=0.46\textwidth]{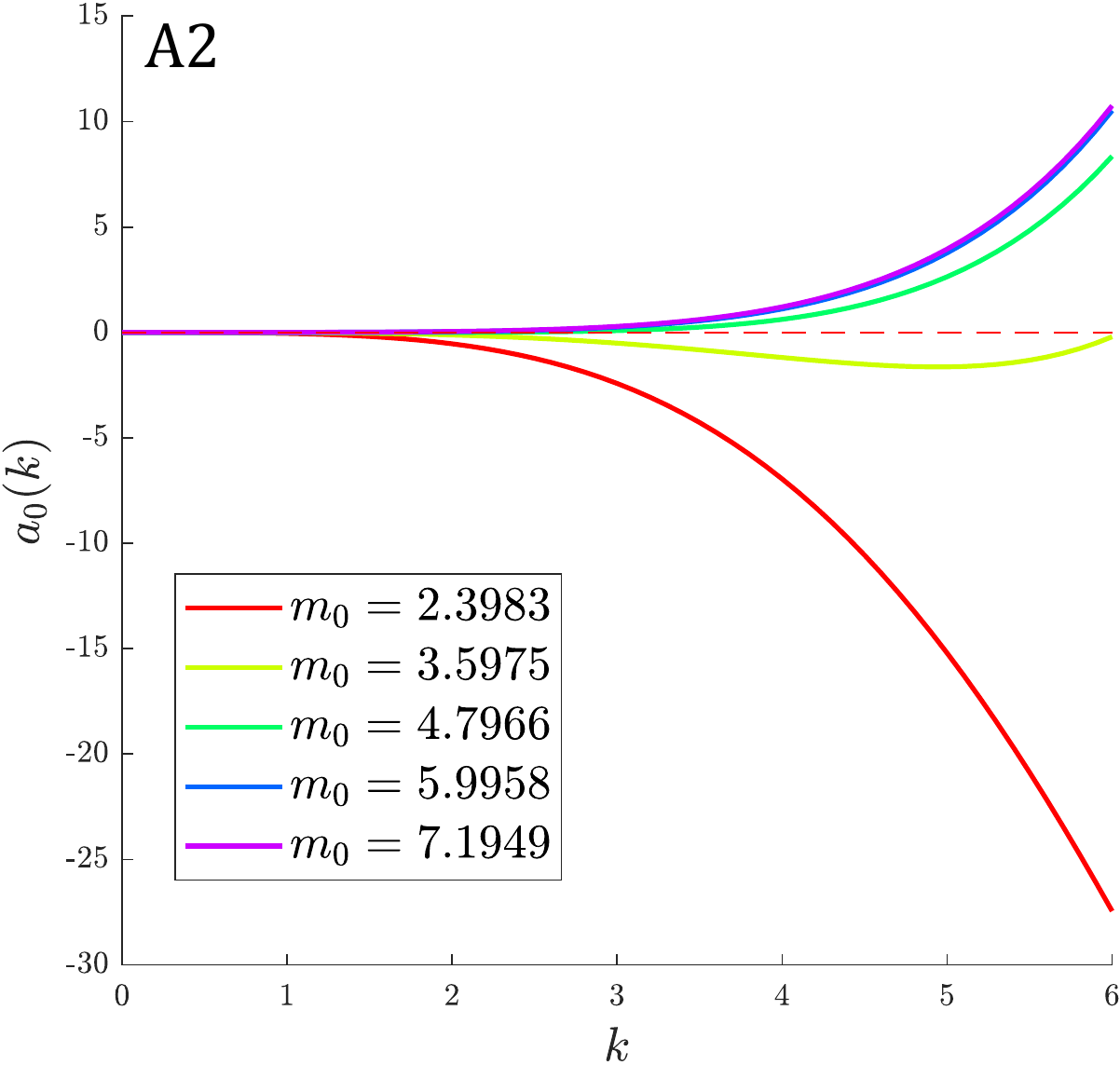}
\includegraphics[width=0.46\textwidth]{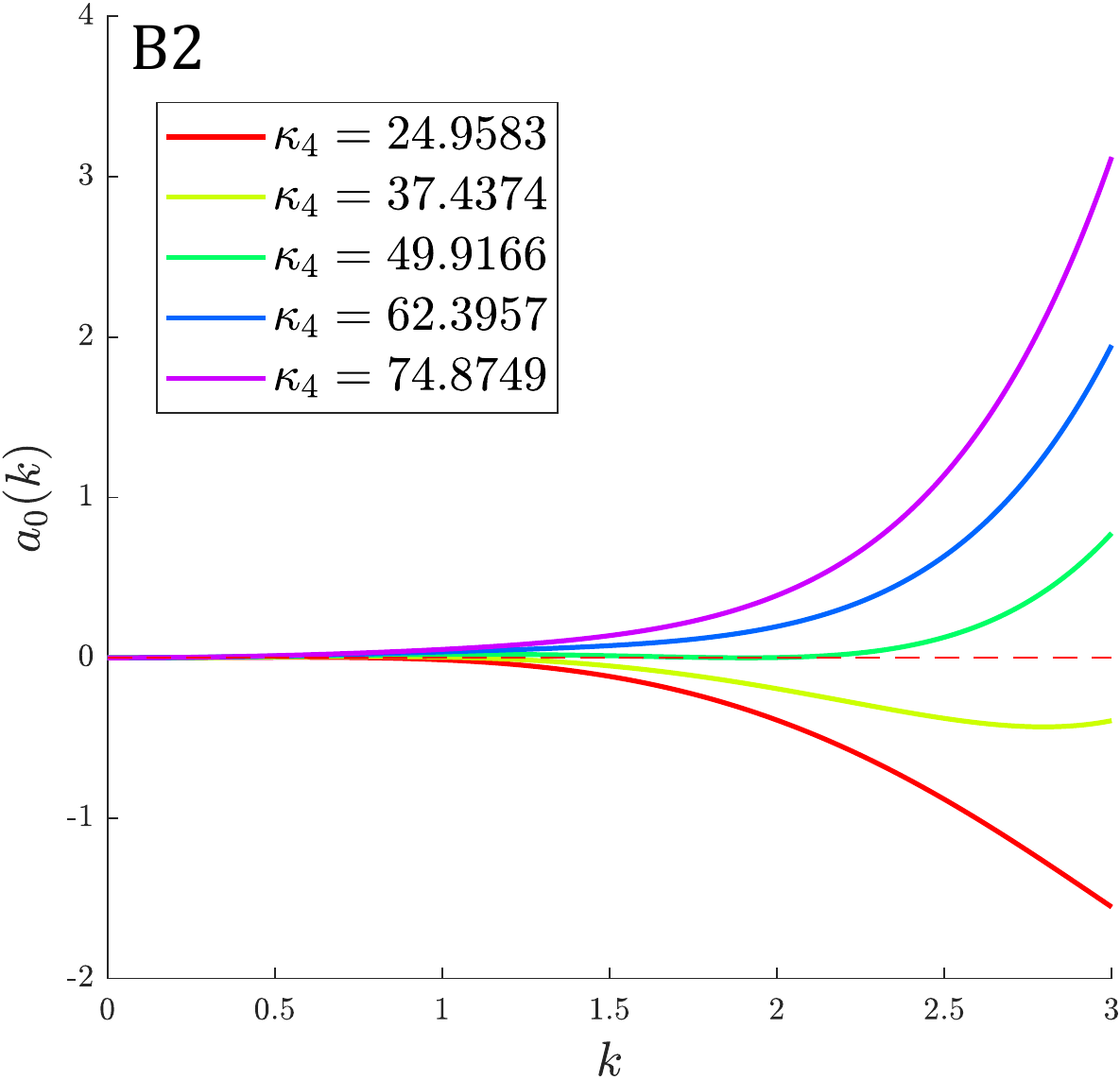}\\
\includegraphics[width=0.46\textwidth]{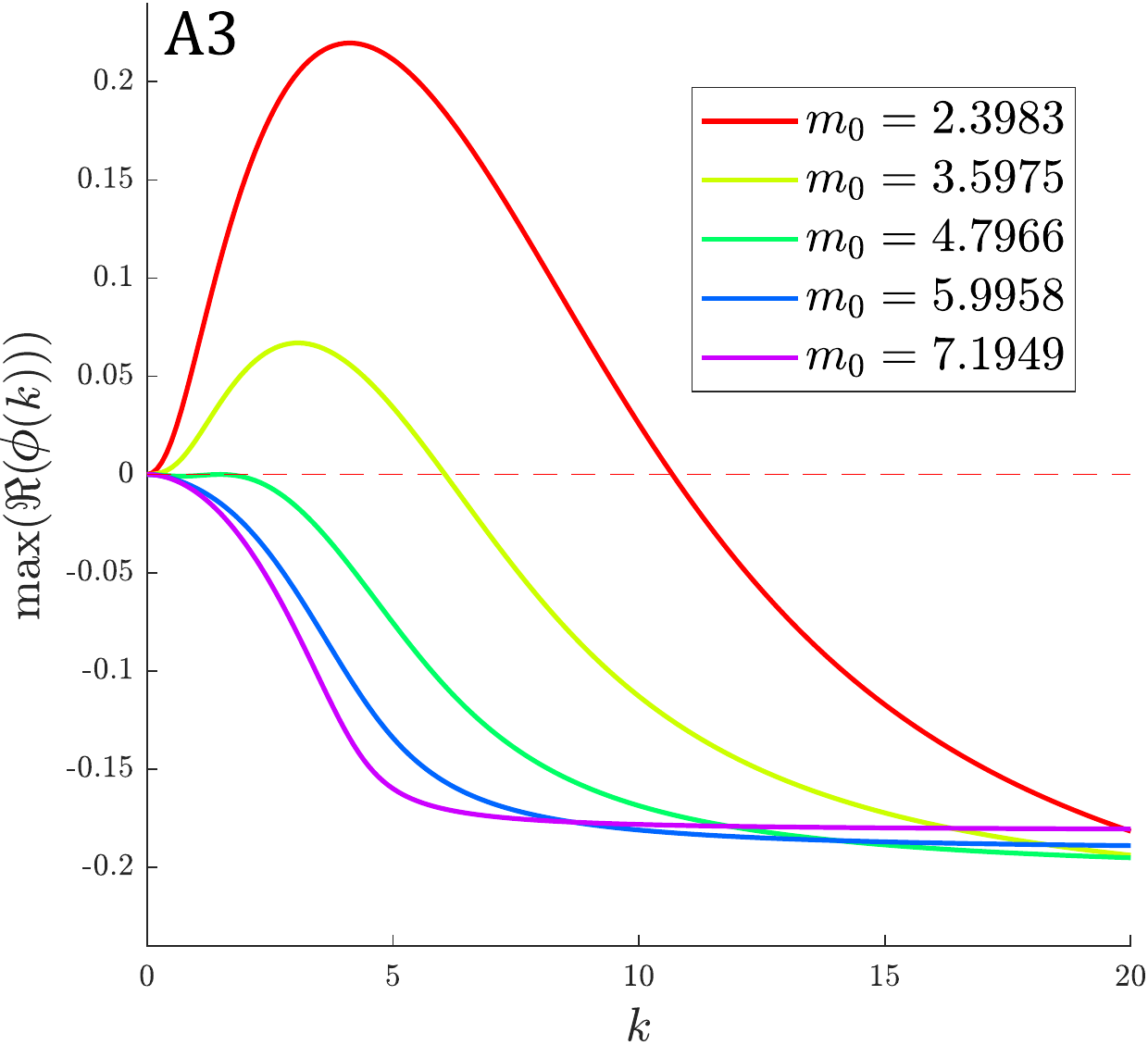}
\includegraphics[width=0.46\textwidth]{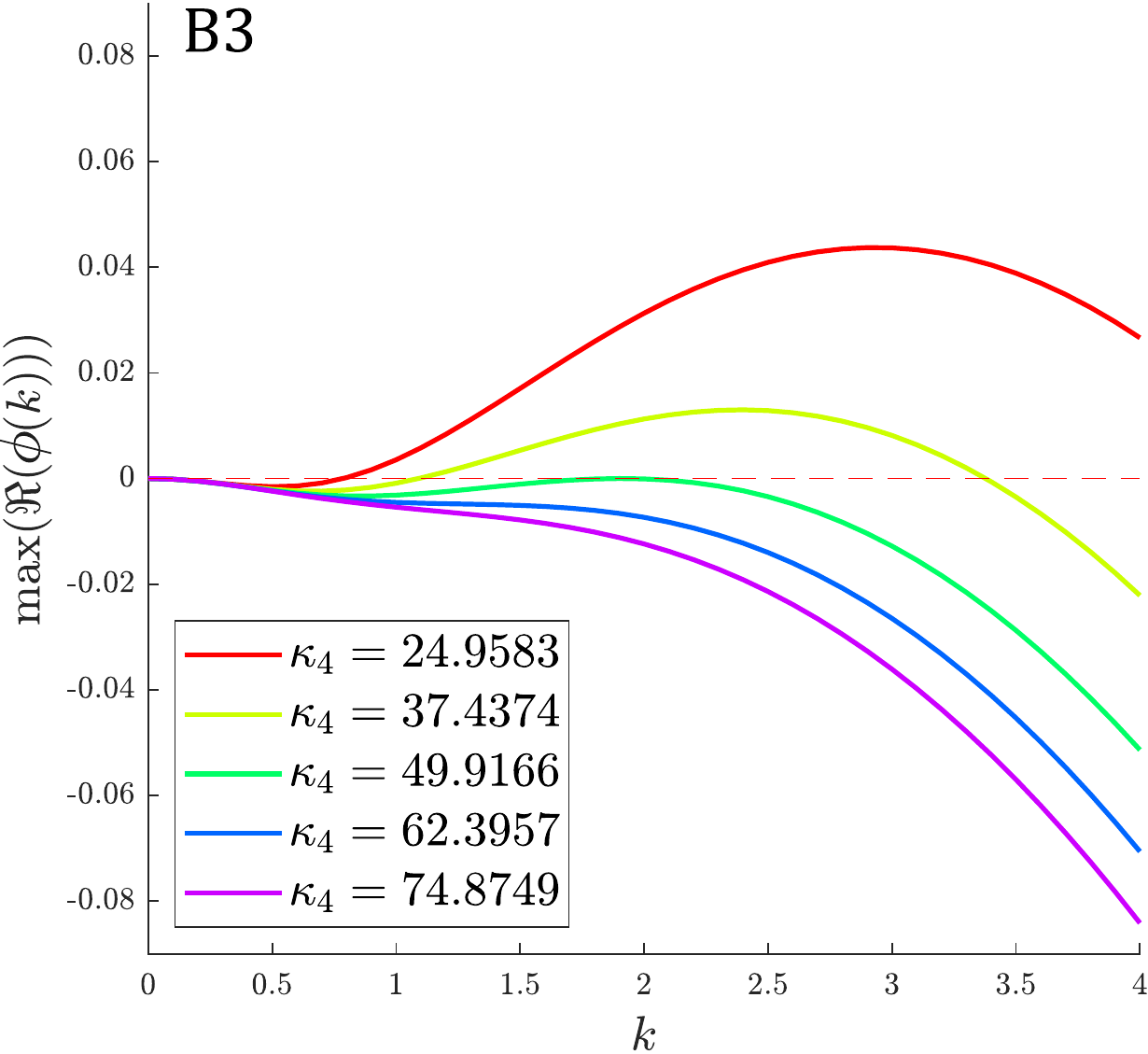}
\end{center}
\caption{Patterning space, parameter condition and dispersion relations for the uncoupled poro-mechano-chemical model. (A1) Predicted pattering space for a selected interval in $(m_0,\alpha)$ parameter space: (green plain) boundary constructed from \eqref{eq:condUC1}; (red-dashed) boundary coming from \eqref{eq:condUC2}; (blue-dot-dashed) boundary built from \eqref{eq:condUC3}. (A2) Parameter coefficient condition $a_0$. Curves are drawn from the critical value $m_0^c$ (green) and for 25\% and 50\% increase/decrease parameter values. (A3) Associated dispersion relations. Colour code is kept identical with (A2). (B1-B3) Similar analysis for the $(\kappa_1,\kappa_4)$ parameter space.}
\label{fig:uncoupledLinearAnalysis}
\end{figure}

Either \eqref{eq:condUC1} or \eqref{eq:condUC2} enforces that $k^2$ has a real positive part, while \eqref{eq:condUC3} is necessary to have a real wave number. These conditions are similar to the ones stated in \cite{painter18}.
The inequality \eqref{eq:condUC2} emphasises how the positive feedback mechanism of the accumulation of mesenchymal cells through chemotaxis leads to stimulate further FGF secretion, based on the activation of the epithelium. When the condition is satisfied, patterning is achieved. 

Condition \eqref{eq:condUC1} indicates that the steady state is destabilised when $A_m>0$, occurring when the clustering-mediated activation is sufficiently strong. Figure \ref{fig:uncoupledLinearAnalysis} (panels (A1)-(B1)) present the patterning space based on the implicit functions defined in \eqref{eq:condUC1}-\eqref{eq:condUC3} for the $(m_0,\alpha)$ and $(\kappa_1,\kappa_4)$ space. The conditions are represented in both Figures by a plain-green curve \eqref{eq:condUC1}, a red-dot-dashed curve \eqref{eq:condUC2}, and a blue-dashed curve \eqref{eq:condUC3}. The boundaries present very similar results to the linear analysis from \cite{painter18}. Here, we accentuate the region by filling with light grey. Panels (A2)-(A3) and (B2)-(B3) in Figure \ref{fig:uncoupledLinearAnalysis} portray the functions $a_0(k^2)$ (parameter condition) and $\lambda(k^2)$ (dispersion relation) for the studied parameters $m_0$ and $\kappa_4$. We present in different colours the behaviour of the system for the critical parameter value (green curve), as well as parameter values obtained by increasing and decreasing by 25\% and 50\% the critical value, where we recall that fixed parameters are taken from Table \ref{tab:LAparams}. We see that moving above or below the critical value results in a finite interval of wave numbers for which the complex eigenvalue presented positive real component and thus leading to instability (patterning). As illustrated also in \cite{painter18}, for a specific parameter set, panels (A1)-(A3) from Figure \ref{fig:uncoupledLinearAnalysis} show that a sufficiently dense dermis cells concentration is necessary to produce instability. While increasing the chemotaxis sensibility of the system we spread the range of possible mesenchymal cells density, by the u-shape form of the boundary, a too-high density cell can lead to prevention of pattern formation (see Figure \ref{fig:uncoupledLinearAnalysis} (A1)). The capacity of BMP to deactivate the epithelium by increasing $\kappa_4$, drops rapidly the patterning ability of the system.

\subsection{Zero activation/inactivation of epithelium.} Different scenarios are possible in such case: only zero activation ($\kappa_1=\kappa_2=0$, $e_0=0$), only zero inactivation ($\kappa_3=\kappa_4=0$, $e_0=1$) or both zero activation and inactivation of the epithelium. Contrary to \cite{painter18}, the coupling to poroelastic structure leads to a more involved analysis of the system and in particular, it does not guarantee that patterning is not predicted. In any case, the coefficients $A_m$ and $A_b$ will be zero and bring the coefficients of $P_2$ to be strictly positive. As the coefficients of $P_3$ are also non-negative, $P(\phi;k^2)$ can only lead to negative eigenvalues if the coefficients of $P_1$ are negative. By Definition \ref{def:Poly}, we know that the $b_i$'s are negative only if the coupled-dependent constant is itself negative. This occurs when the following conditions are satisfied
\begin{linenomath*}
\[
-\taup > 0 \Leftrightarrow \Bigg\lbrace \begin{matrix}
\tau > 0 & \text{ and } &  m_0 < \sqrt{\frac{1}{\zeta}} \\
\tau < 0 & \text{ and } &  m_0 > \sqrt{\frac{1}{\zeta}}
\end{matrix}.
\]
\end{linenomath*}
If one of the condition is satisfied, the formation of patterns is expected. This is an opposite conclusion as that drawn in the case of chemotaxis only \cite{painter18}.


\begin{figure}[!t]
\begin{center}
\includegraphics[width=0.325\textwidth]{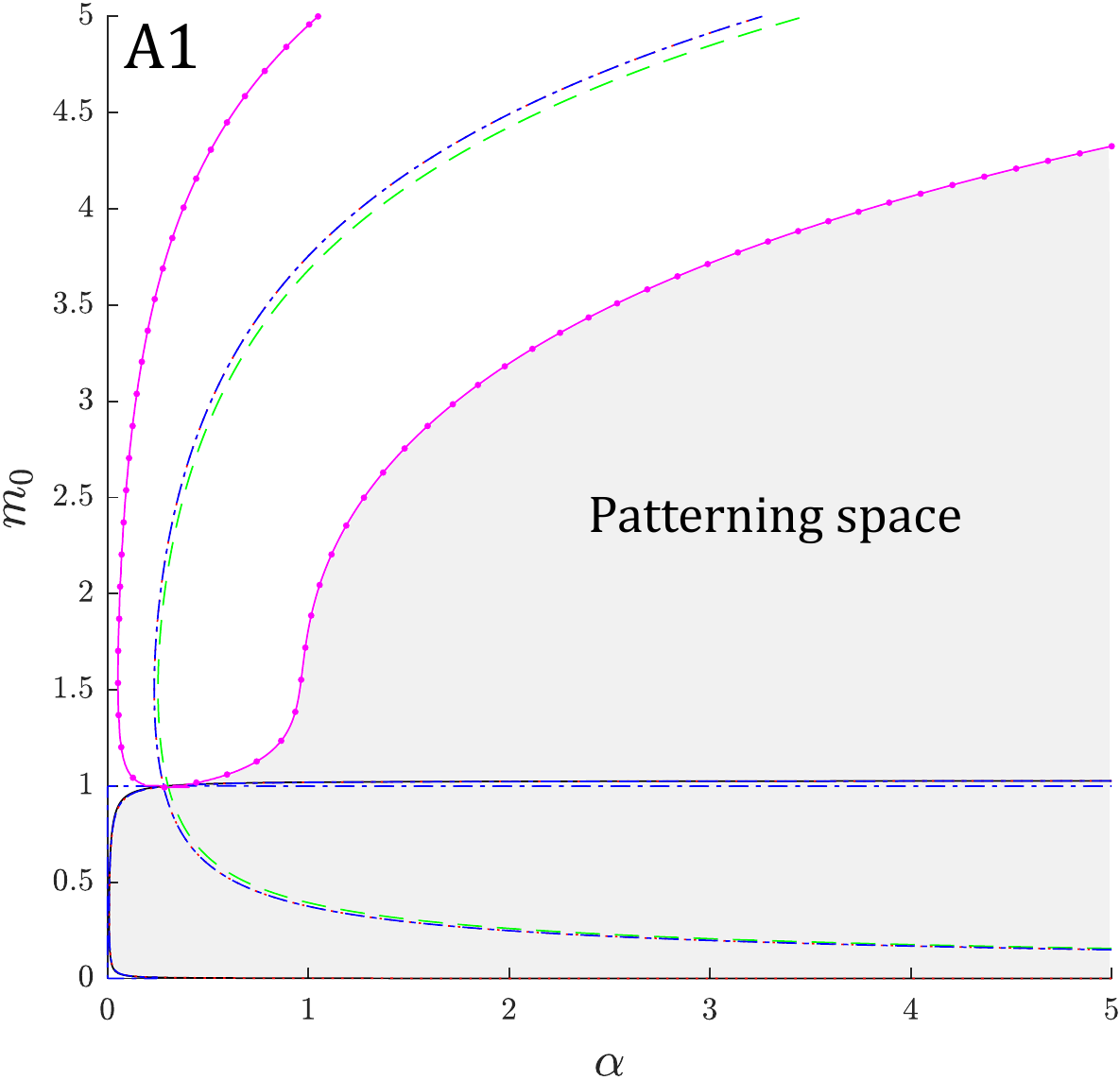}
\includegraphics[width=0.325\textwidth]{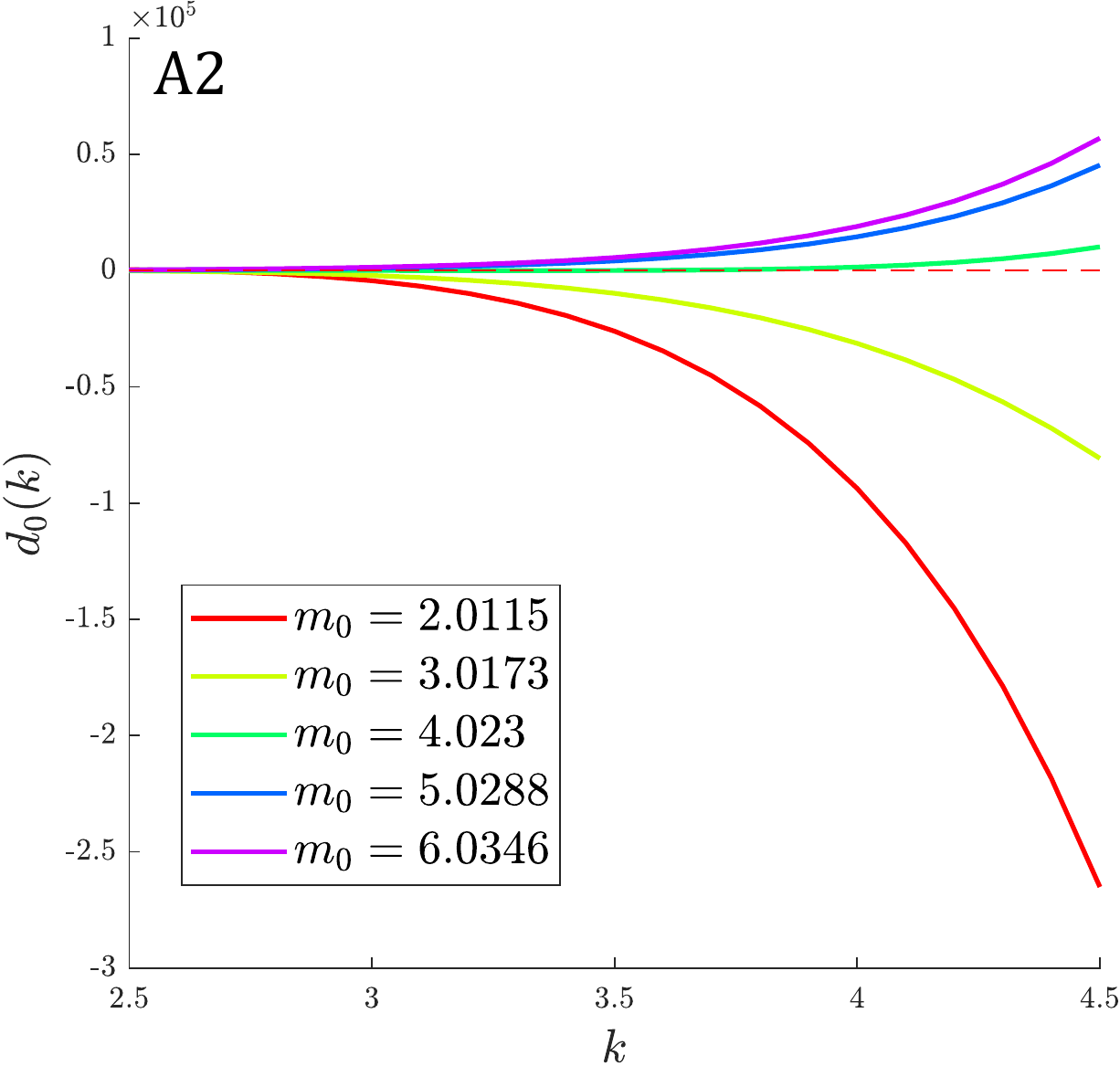}
\includegraphics[width=0.325\textwidth]{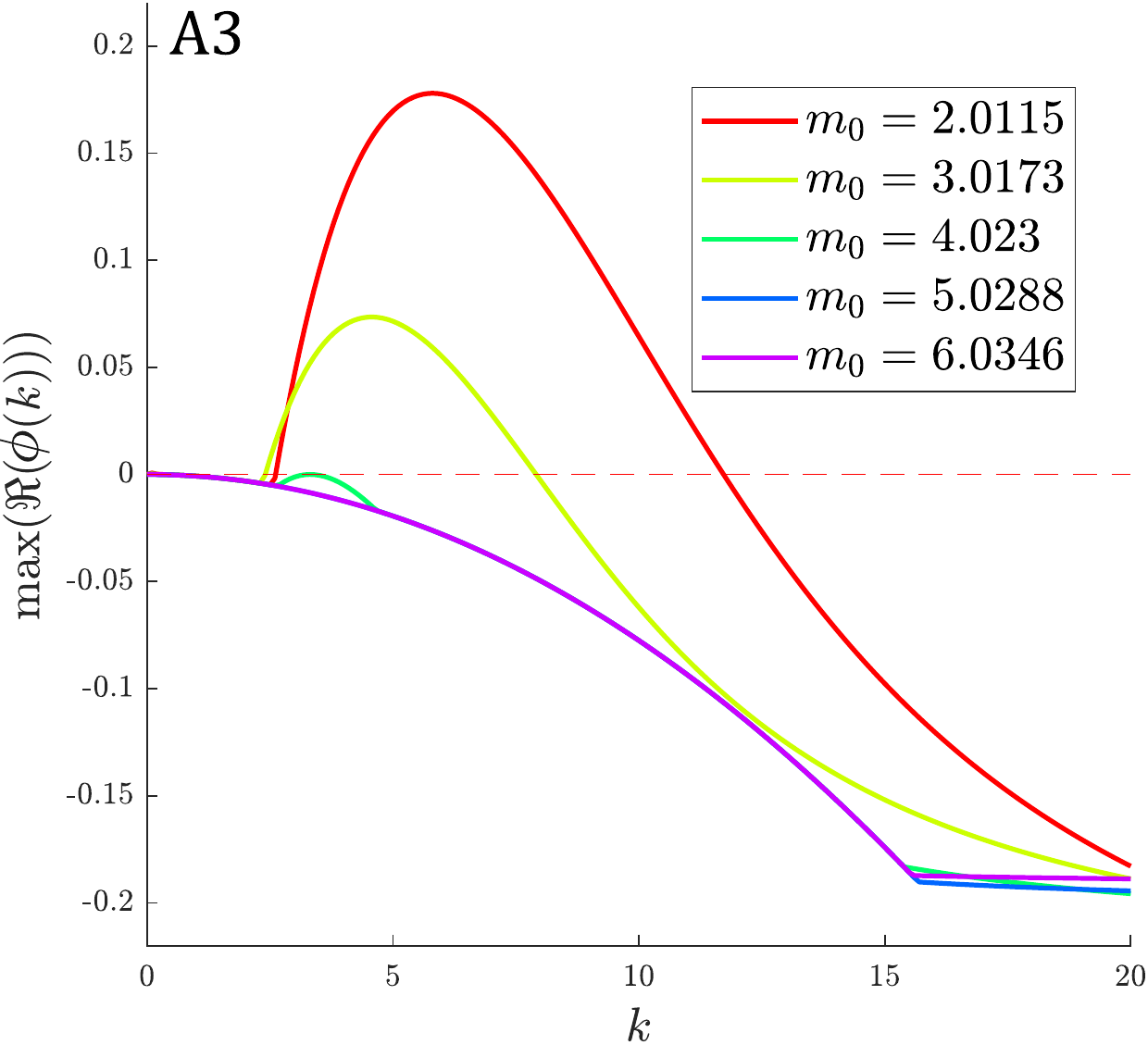}\\
\includegraphics[width=0.325\textwidth]{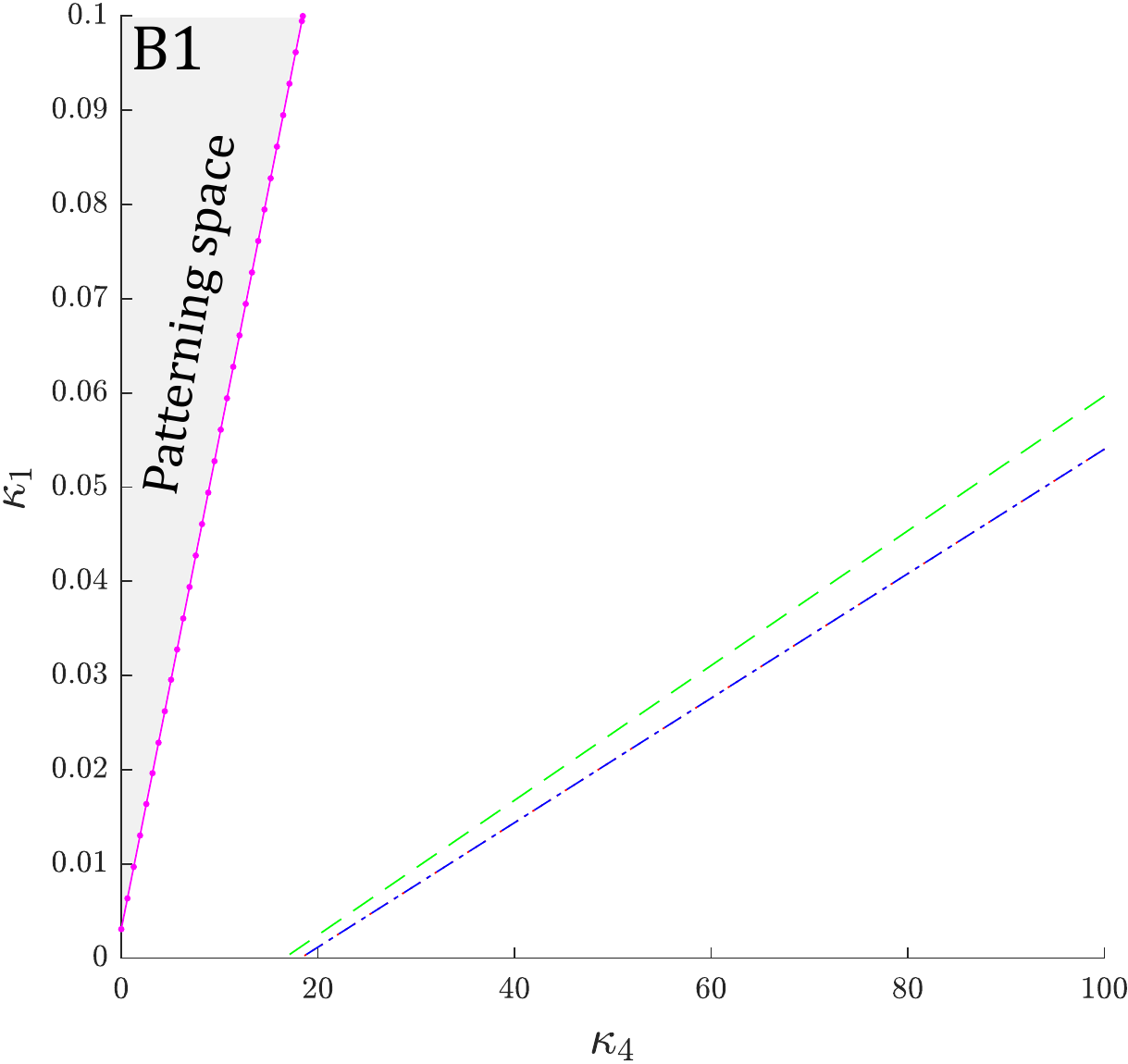}
\includegraphics[width=0.325\textwidth]{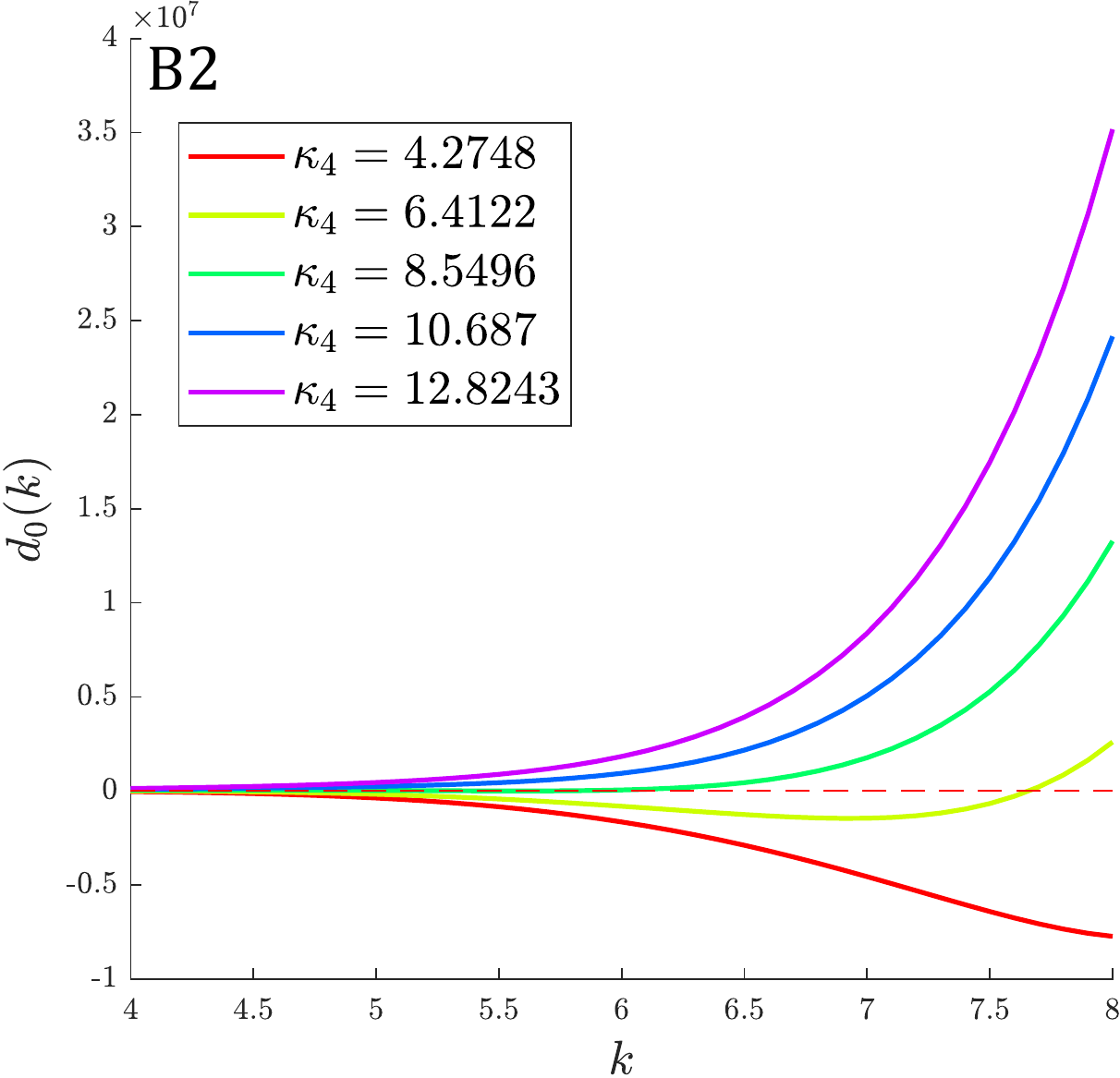}
\includegraphics[width=0.325\textwidth]{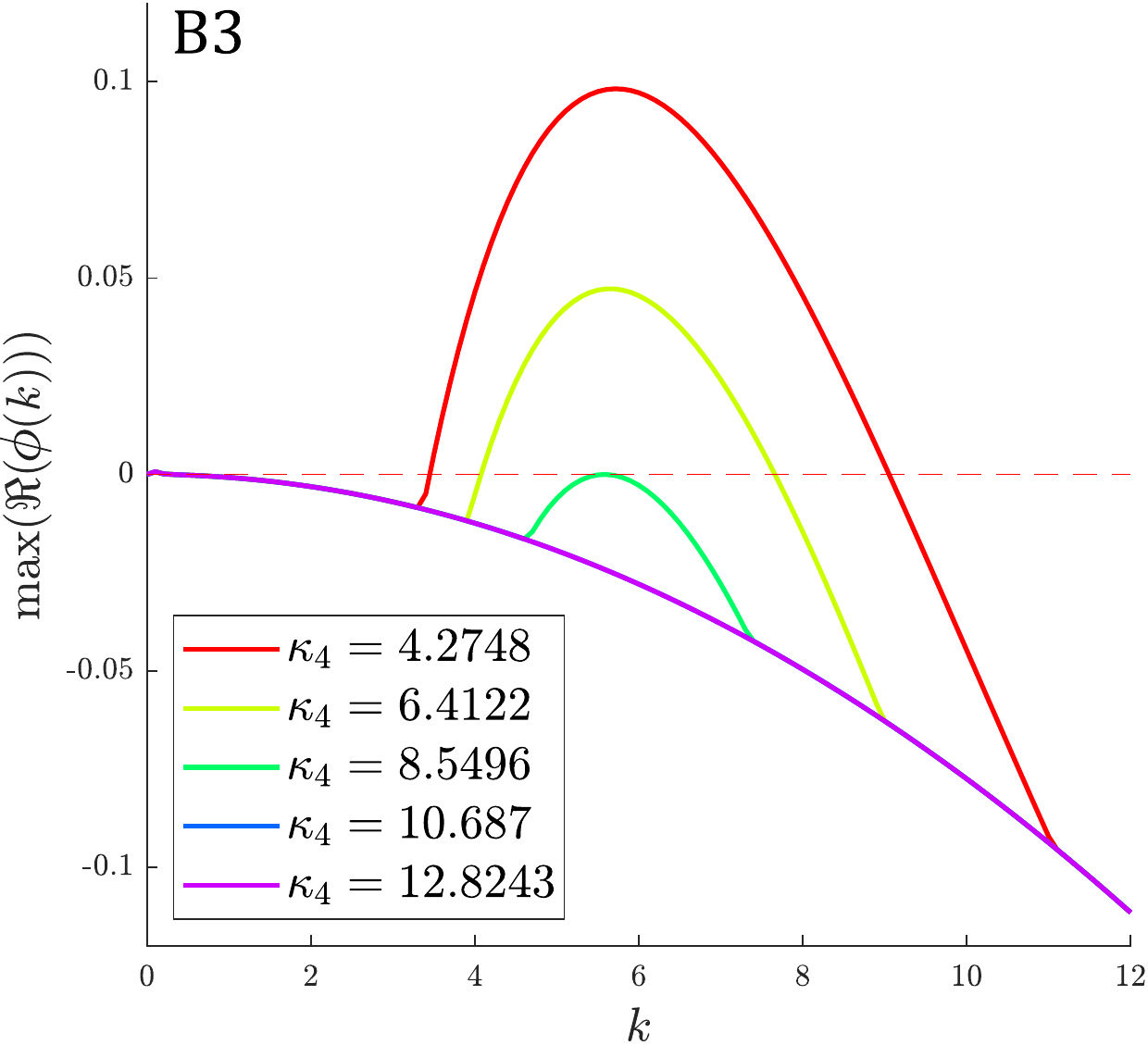}\\
\includegraphics[width=0.325\textwidth]{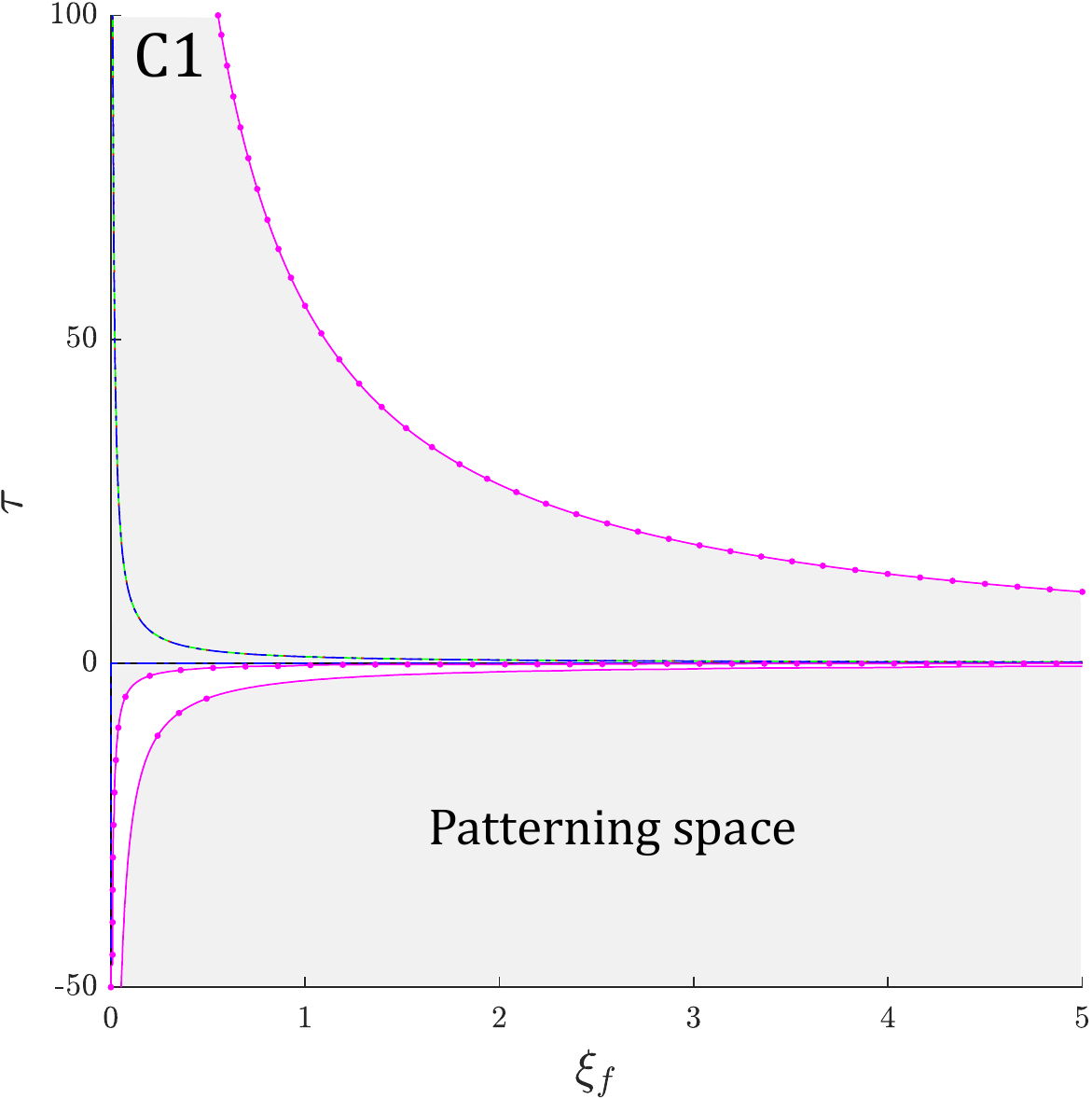}
\includegraphics[width=0.325\textwidth]{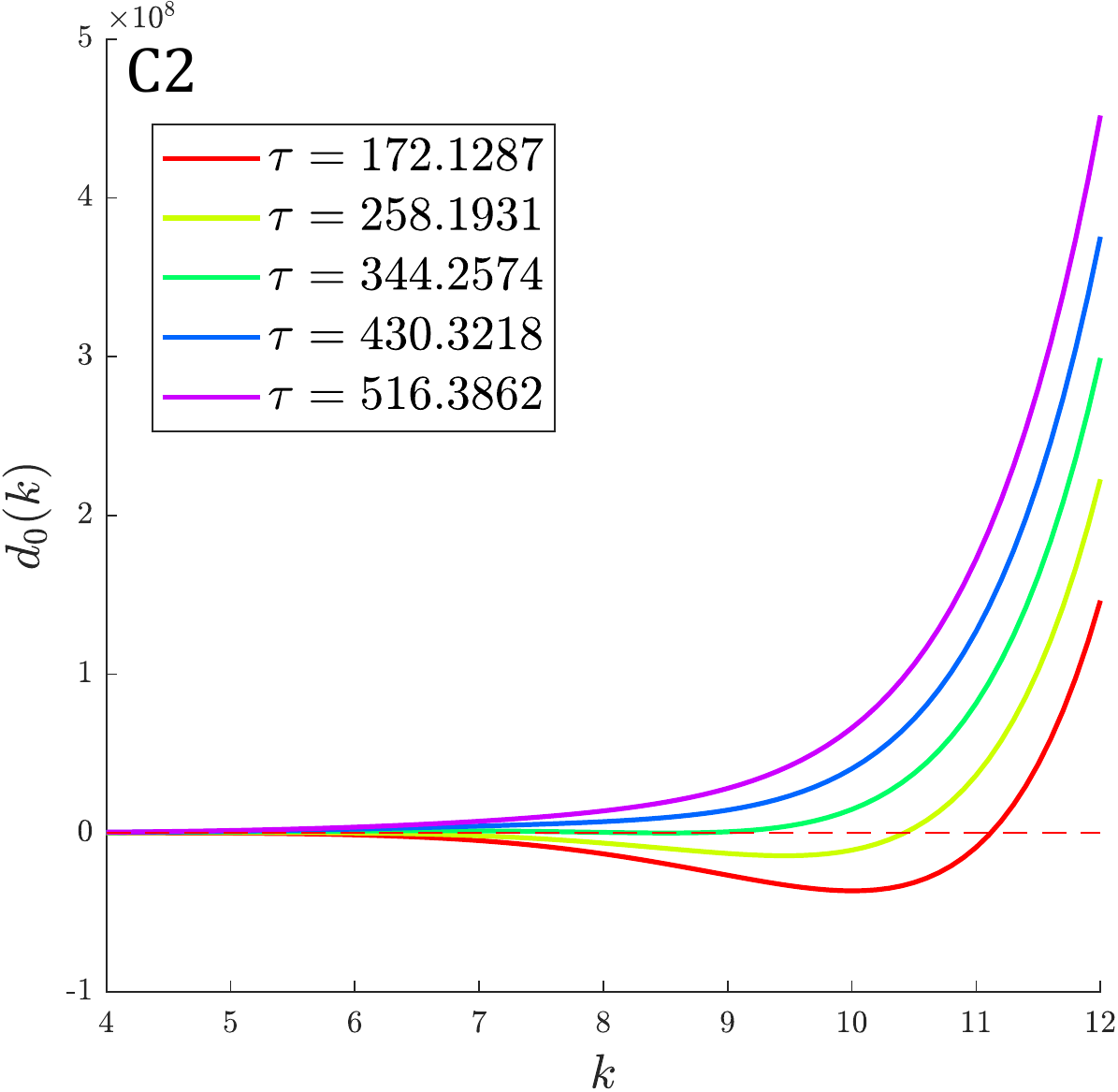}
\includegraphics[width=0.325\textwidth]{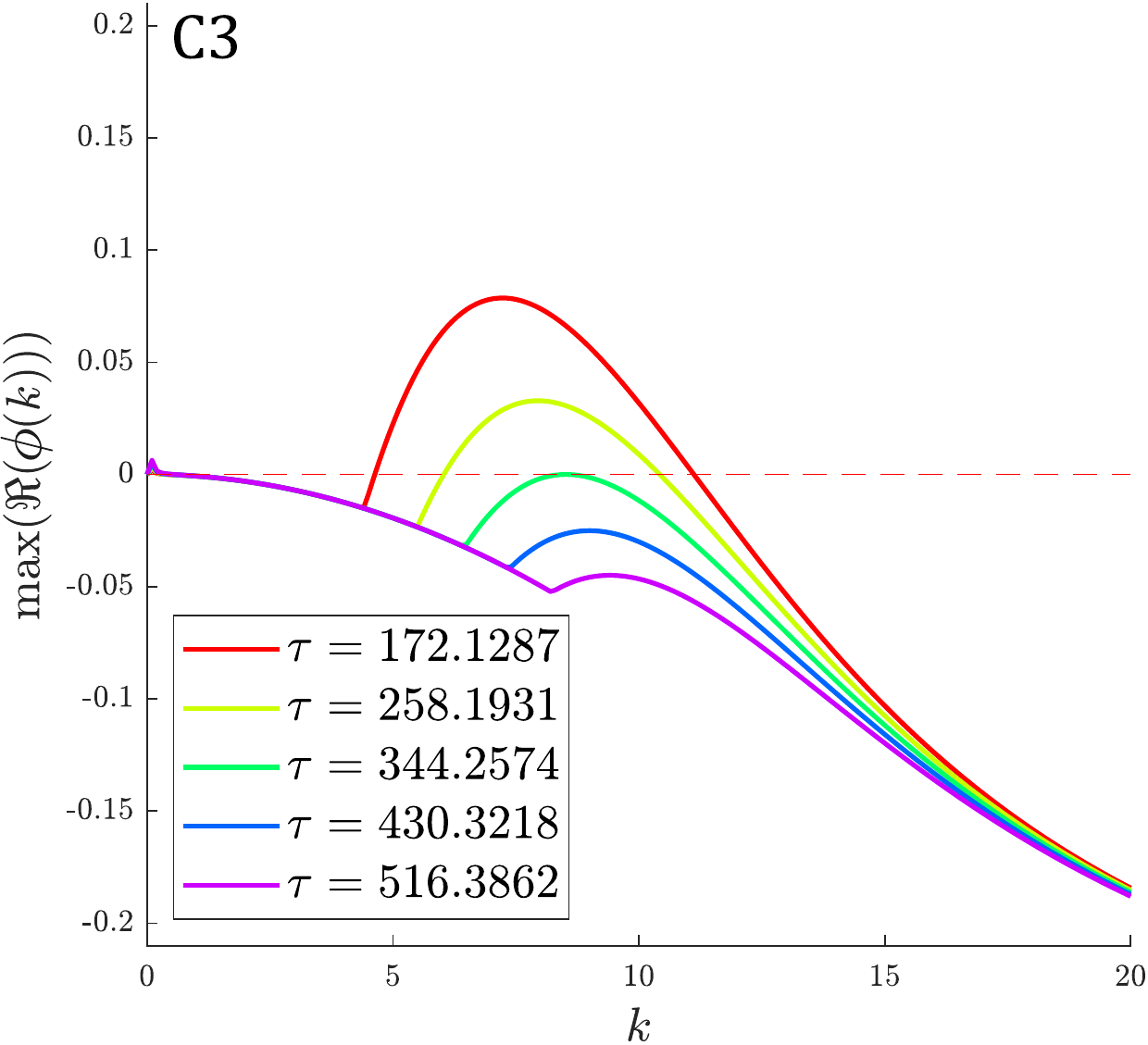}
\end{center}
\caption{Patterning space, parameter condition and dispersion relations for the coupled poro-mechano-chemical model. (A1):  predicted pattering space for a selected interval in $(m_0,\alpha)$. Parameter space: (black plain) boundary built from \eqref{eq:condC1} for $\theta_1$; (green-dashed) boundary built from \eqref{eq:condC1} for $\theta_2$; (red-dot-dot) boundary built from \eqref{eq:condC2}; (blue-dot-dashed) boundary from \eqref{eq:condC3}; (magenta-dot-plain) boundary built from \eqref{eq:condCDisc}. (A2): Coefficient condition on $d_0$. Curves are drawn from the critical value $m_0^c$ (green) and for the 25\% and 50\% increase/decrease parameter value. (A3): associated dispersion relations. Colour code is kept identical with (A2). (B1-B4): similar analysis for the $(\kappa_1,\kappa_4)$ parameter space. (C1-C3): similar analysis for the $(\tau,\xi_f)$ parameter space.}
\label{fig:coupledLinearAnalysis}
\end{figure}

\subsection{General case.} The thorough analysis of the coupled system involves a seventh-order polynomial. Analytical Routh--Hurwitz conditions in this case are hardly accessible and the crossed effects due to multiple parameters and conditions will inevitably make difficult the interpretation of the system behaviour. We can however restrict then the analysis to $\phi^0$, with coefficients $d_0=a_0c_0+b_0$, where $a_0$, $b_0$, $c_0$ are as in Definition \ref{def:Poly}. After some algebraic manipulations, we derive $d_0$ with respect to $k^2$ and obtain
\begin{linenomath*}
\[
d_0^{\prime}(k^2) = \theta_4 k^{8} + \theta_3 k^6 + \theta_2 k^4 + \theta_1 k^2,
\]
\end{linenomath*}
where
\begin{linenomath*}
\begin{align*}
\theta_1	 &= - 2 \xi_f \taup \Ep ( \kon + \koff ),\\
\theta_2 &= 3 \Bigg[ ( 2\mu + \lambda )\Big( D_m \delta_B ( \kon + \koff ) + \Ep H_3 A_b - \Ep A_m \Big)  \\
& \qquad - \xi_f \taup \Ep ( \kon + \koff ) \Bigg],\\
\theta_3 &= 4 ( 2\mu + \lambda )\left( (D_m \delta_B + D_m D_f )( \kon + \koff ) - \Ep A_m \right),\\
\theta_4 &= 5 D_m D_f( 2\mu + \lambda )( \kon + \koff).
\end{align*}
\end{linenomath*}

\begin{figure}[!t]
\begin{center}
\includegraphics[width=0.5\textwidth]{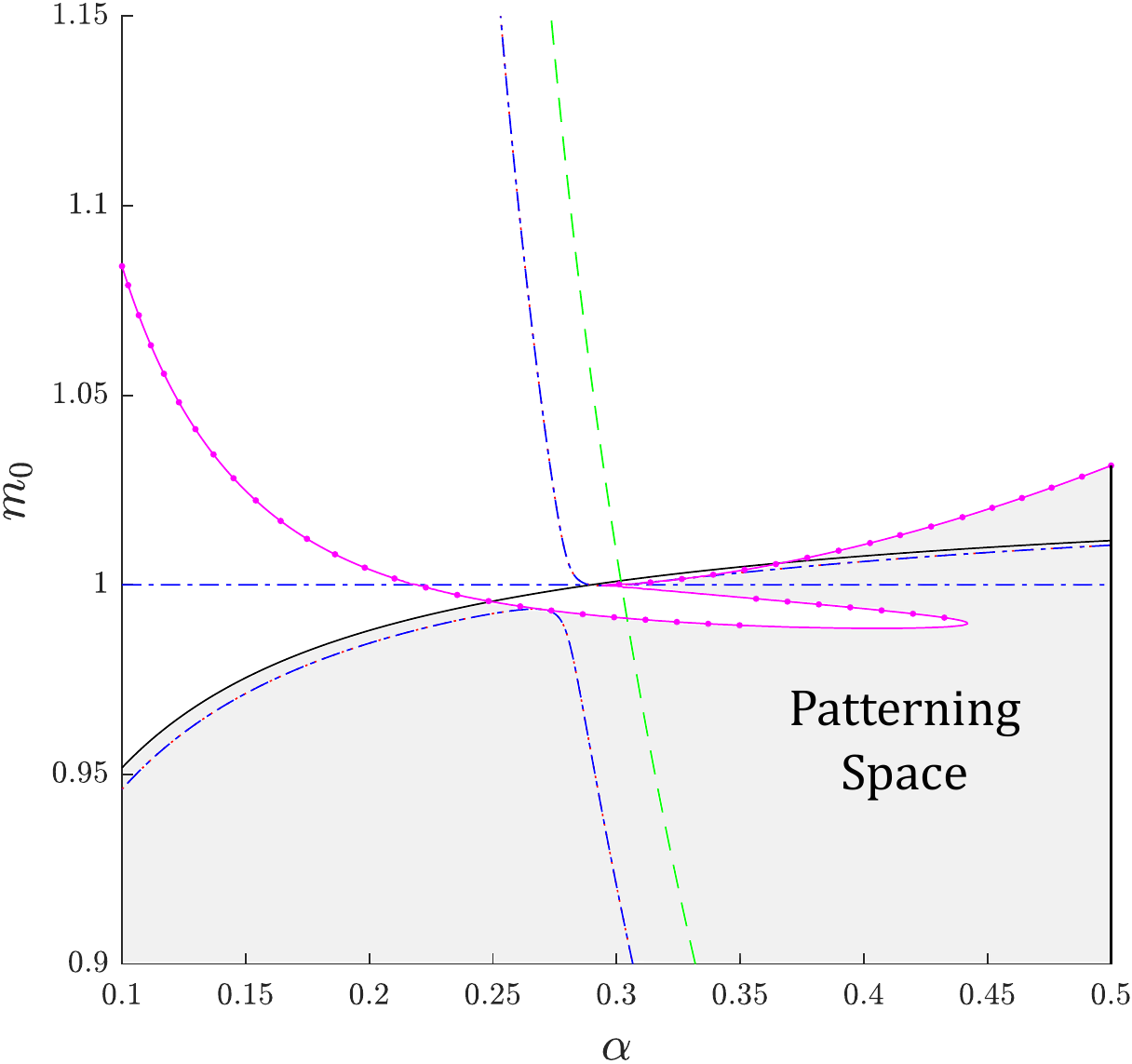}
\caption{Patterning space for the coupled poro-mechano-chemical model. Selected interval in the $(m_0,\alpha)$ parameter space. (black plain) boundary built from \eqref{eq:condC1} for $\theta_1$; (green-dashed) boundary from \eqref{eq:condC1} for $\theta_2$; (red-dot-dot) boundary from \eqref{eq:condC2}; (blue-dot-dashed) boundary from \eqref{eq:condC3}; (magenta-dot-plain) boundary built from \eqref{eq:condCDisc}.\label{fig:patternSpaceZoom}}
\end{center}
\end{figure}

In contrast with the uncoupled case studied above, here $k_c^2$ is found by solving a quartic polynomial (employing the intrinsic MATLAB function \texttt{roots}, after constructing the associated vector of polynomial coefficients). In addition, one needs to satisfy some additional conditions leading to $d_0<0$ for any positive $k^2$. More specifically, the constrains consist on at least one of the Routh--Hurwitz conditions \cite{routh1877}, stated as
\begin{linenomath*}
\begin{subequations}
\begin{align}
	\theta_1 > 0, \qquad \theta_2 > 0, \qquad \theta_3 & > 0, \label{eq:condC1} \\
	\theta_2\theta_3 - \theta_1\theta_4 & > 0, \label{eq:condC2} \\
\theta_1\theta_2\theta_3 - \theta_1^2\theta_4 & > 0, \label{eq:condC3}
\end{align}\end{subequations}
\end{linenomath*}
or the discriminant
\begin{linenomath*}
\begin{equation}
	 -27\theta_4^2\theta_1^4 + 18\theta_4\theta_3\theta_2\theta_1^3 - 4\theta_4\theta_2^3\theta_1^3 - 4\theta_3^3\theta_1^3 + \theta_3^2\theta_2^2\theta_1^2 > 0 \label{eq:condCDisc}.
\end{equation}
\end{linenomath*}
Condition \eqref{eq:condC2} highlights the effect of poromechanics towards the instability of the system, in comparison against the uncoupled cased. Increasing the influence of the mesenchymal cells through $\tau$, or the influence of mechanics on $f$ through $\xi_f$, can relax the requirement that the feedback-positive loop from \eqref{eq:condUC2} should be sufficiently strong to create instability. As mentioned in the activation/inactivation of epithelium scenario, this will depend mostly on the sign of the active stress component of the coefficient, and thus on the sign of $\tau$ and the value of $m_0$, see Figure~\ref{fig:coupledLinearAnalysis}(C1). The dispersion relation for the $(m_0,\alpha)$- and $(\kappa_1,\kappa_2)$-- parameter spaces illustrates how, for the same fixed set of parameters from Table \ref{tab:LAparams}, the poromechanical coupling leads to a decrease of  the critical parameter value (Figures  \ref{fig:coupledLinearAnalysis}(A3) and (B3)). Unfortunately, the ease to push the system to instability is counter-balanced by a more restrictive space of parameters values where patterns can be formed (more precisely, compared to Figures \ref{fig:uncoupledLinearAnalysis}(A1),(B1) and \ref{fig:coupledLinearAnalysis} (A1),(B1),(C1)). This can increase the difficulty in tuning parameters to produce specific patterns. Indeed, zooming into the parameters-pair $(m_0,\alpha)$ on the region where the boundary conditions cross (see Figure \ref{fig:patternSpaceZoom}), exposes how easily the patterns can disappear after a relatively small variation in parameter values.

\section{Extension to finite-strain poroelasticity and growth}\label{sec:growth}

We briefly recall some kinematic considerations, needed in order to incorporate growth effects. Thorough
reviews of mechanical models for growth of soft living tissues can be found in \cite{jones12,kuhl14}.
Let us
denote by $\bX$ the position of a material point in the reference (or initial) configuration (region  $\Omega_0$ with boundary $\partial\Omega_0 = \Gamma_0 \cup \Sigma_0$ and outward pointing unit normal $\bN_0$) 
and let $\bx=\boldsymbol{\varphi}(\bX)$ denote its position at time $t$ in the current configuration
of the domain, $\Omega$, determined by the deformation mapping $\boldsymbol{\varphi}(\cdot,t): \Omega_0\to \Omega$.  We adopt the notation $\bDiv,\vDiv$ for the divergence of tensor and vector fields, 
and $\bGrad, \vGrad$ for the gradients of vector and scalar fields, respectively; all with respect to the material coordinates. The
geometric deformation tensor is $\bF = \bGrad\boldsymbol{\varphi} = \bI + \bGrad\bu$ (where the gradient of the deformation
map and of the displacement field is
taken with respect to the reference coordinates) and its Jacobian is $J = \det \bF > 0$. The velocity in Lagrangian coordinates is denoted as $\bv = \dot{\bu}$, and $\bC= \bF^\intercal\bF$
is the right Cauchy--Green strain tensor. 

As in the linearised regime, we assume that the region occupied by the tissue is fully saturated, and consisting of a fluid and of a solid phase: 
\begin{linenomath*}
\begin{equation}\label{p:phases}
\phi_f+\phi_s=1,\end{equation}
\end{linenomath*}
where the $\phi_j$'s denote the volume fraction of each phase. The phase $s$ is assumed a neo--Hookean solid whereas the liquid phase is considered as an ideal fluid. Both phases are assumed intrinsically incompressible. The expected growth is due to proliferation of cells, and implying a mass exchange between the phases. The mass balances in the Eulerian framework assume the form
\begin{linenomath*}
\begin{subequations}\begin{align}
\label{p:mass1}    \frac{\partial \phi_s}{\partial t} + \vdiv (\phi_s\bv_s) &= r_s,\\
\label{p:mass2}    \frac{\partial \phi_f}{\partial t} + \vdiv (\phi_f\bv_f) &= -r_s,
\end{align}\end{subequations}\end{linenomath*}
with $\bv_s,\bv_f$ the Eulerian velocities in each phase, and $r_s=\alpha_0\phi_s(1-\phi_s)(m-m_0)$ being the mass growth rate (net proliferation rate of cells in the solid phase), here assumed linearly dependent on the concentration of mesynchemal cells and decreasing as the available space decreases. In turn, the general form of the momentum balances is 
\begin{linenomath*}\begin{align*}
    \bdiv (\phi_s \bsigma_{\mathrm{eff},s}) - \phi_s\nabla p + \ff_{sf} & = \cero, \\
    \bdiv (\phi_f \bsigma_{\mathrm{eff},f}) - \phi_f\nabla p + \ff_{fs} & = \cero,
\end{align*}\end{linenomath*}
where the $\bsigma_{\mathrm{eff},j}$'s are the respective effective Cauchy stresses, $p$ is the average pressure of the liquid phase, and
\begin{linenomath*}
  \[\ff_{sf} = p \nabla \phi_s + \eta/\kappa \phi_s \phi_f (\bv_f - \bv_s),\quad \ff_{fs} =p \nabla \phi_f - \eta/\kappa \phi_f \phi_s (\bv_f - \bv_s),\]\end{linenomath*}
are the respective  drag terms and net sources of momentum, where $\kappa,\eta$ are the permeability and dynamic viscosity, respectively. The effective stress of the liquid phase $\bsigma_{\mathrm{eff},f}$ is assumed negligible with respect to the interstitial pressure gradient and therefore the mixture stress is $\bsigma_{\mathrm{eff}}=\bsigma_{\mathrm{eff},s}$. In addition, Darcy's law gives
\begin{linenomath*}
\begin{equation}\label{p:Darcy}
\bv_f-\bv_s = - \frac{\kappa}{\eta\phi_f} \nabla p,
\end{equation}
  \end{linenomath*}
and then the mixture momentum balance is written as
\begin{linenomath*}
\begin{equation}\label{p:momentum}- \bdiv(\bsigma_{\mathrm{eff}}) + \nabla p = \cero \qquad \text{in }\Omega\times(0,t_{\mathrm{final}}].\end{equation}
    \end{linenomath*}
    
\begin{figure}[!t]
\begin{center}
\includegraphics[width=0.65\textwidth]{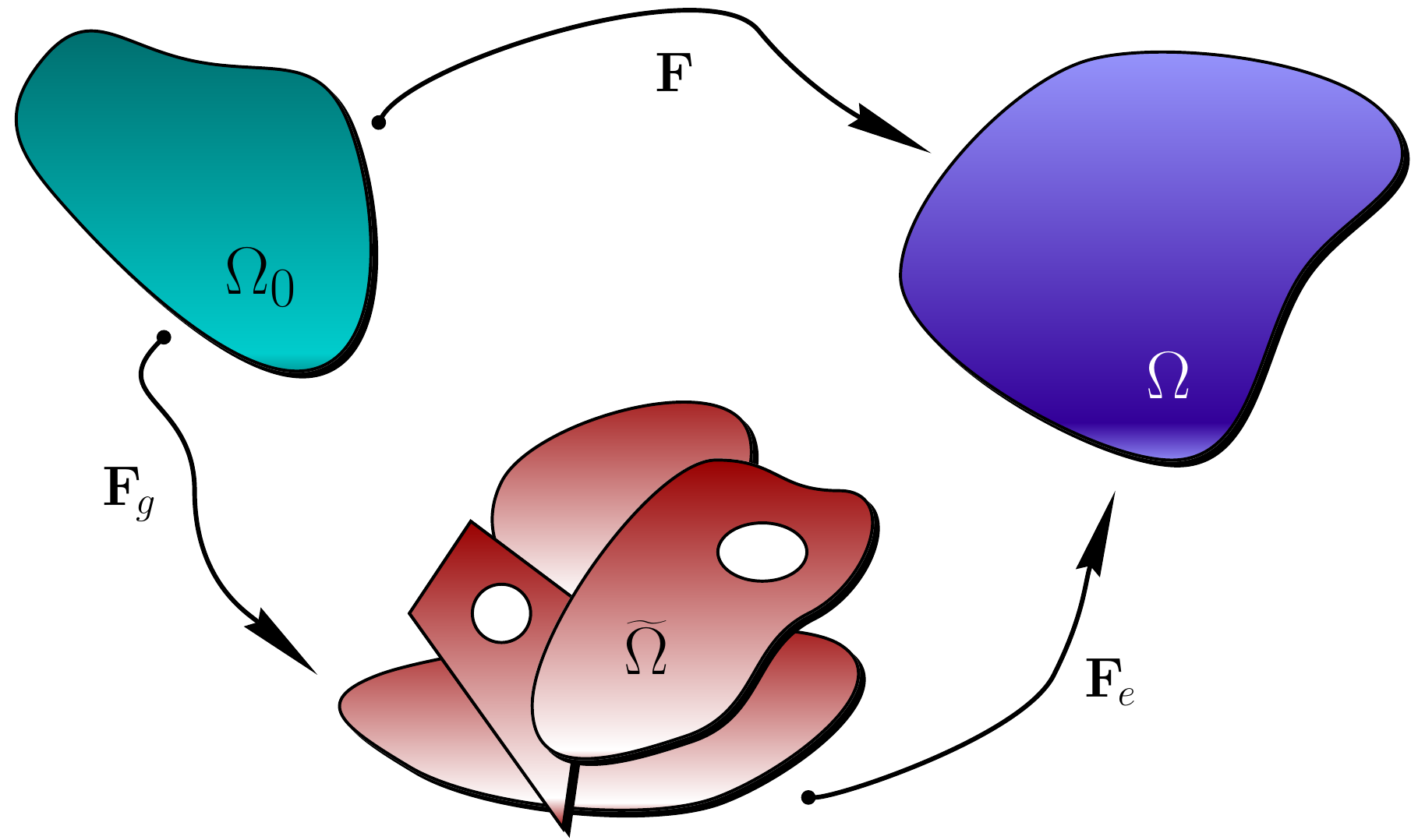}
\caption{Decomposition of the tensor gradient of deformation $\bF$ into pure growth $\bF_g$ and an elastic deformation tensor $\bF_e$.
The intermediate configuration $\widetilde{\Omega}$, between the undeformed reference state $\Omega_0$ and
the current/final configuration $\Omega$ including growth and elastic response with stress, is an incompatible and stress-free growth state.\label{fig:decomposition}}
\end{center}
\end{figure}

The deformation mapping can be decomposed into a purely growth part and the remainder, elastic, deformation.
Such splitting implies that there exists an intermediate configuration $\widetilde{\Omega}$ between
$\Omega_0$ and $\Omega$ which is not necessarily compatible (and which we assume is completely stress free, see sketch in
Figure~\ref{fig:decomposition}),
and consequently a multiplicative decomposition
of the  deformation gradient into a growth deformation gradient (describing local generation or removal of material points)
and an elastic deformation gradient tensor is admitted
(see, for instance, \cite{lee69,rodriguez94,kida18})
\begin{linenomath*}
\begin{equation}\label{eq:FeFg}
\bF = \bF_e\bF_g,
\end{equation}
\end{linenomath*}
and therefore $J = J_eJ_g$ with $J_e := \det \bF_e$, $J_g := \det \bF_g$. Note also that since $\bF$ is assumed nonsingular, so are the tensors $\bF_e,\bF_g$.

Since the intermediate configuration is considered stress-free, stresses are exerted only by
the elastic deformation, and the constitutive relations between a given strain energy function $\Psi$
(that characterises the material response of the solid for hyperelastic materials)
and the measures of stress
can be stated
with respect to the intermediate configuration  $\widetilde{\Omega}$. For the Cauchy stress this gives
\begin{linenomath*}
\begin{equation}\label{p:stress}
\bsigma_{\mathrm{eff}}= 2J_e^{-1}\bF_e \frac{\partial\Psi}{\partial\bC_e}\bF_e^{- \intercal},\end{equation}
  \end{linenomath*}
  where $\bC_e = \bF_e^\intercal\bF_e$ and 
$J_e = \det (\bF\bF_g^{-1})$ is the reversible part of the Jacobian. This constitutive equation requires to describe also $\bF$ and $\bF_g$, but from the relation $\dot{\bF}{\bF^{-1}} = \bnabla\bv_s$, we realise that only $\bF_g$ needs to be specified.

In order to translate equations \eqref{p:mass1}-\eqref{p:stress} into Lagrangian form, we recall the change in volume and area relations $\mathrm{d}V = J\mathrm{d}V_0$ and $\mathrm{d}S = J\bF^{-\intercal}\mathrm{d}S_0$ together with the transformation of unit outward normal vectors $\nn = \bF^{-\intercal}\bN_0$. For the solid mass balance \eqref{p:mass1} it suffices to apply Reynolds transport theorem, the change of coordinates, and Maxwell's localisation to obtain
\begin{linenomath*}
\begin{equation}\label{q:mass1}
  \dot{\overline{J\phi_s}} = J r_s.\end{equation}
  \end{linenomath*}
The Lagrangian form of Darcy's law   
\begin{linenomath*}
\begin{equation}\label{q:Darcy}
\bv_f-\bv_s = - \bF\frac{\kappa}{\eta\phi_f}\bC^{-1}\vGrad p,
\end{equation}
  \end{linenomath*}
results from \eqref{p:Darcy} in combination with transport, localisation, as well as pull-back operations.   
For the liquid mass balance \eqref{p:mass2}, one can apply the generalised Reynolds transport theorem, \eqref{q:Darcy},  
divergence's theorem, and localisation, to arrive at 
\begin{linenomath*}
\begin{equation}\label{q:mass2}
 \dot{\overline{J\phi_f}} + \vDiv(J\frac{\kappa}{\eta}\bC^{-1}\vGrad p) = -J r_s.\end{equation}
  \end{linenomath*}
Summing up \eqref{q:mass1} and \eqref{q:mass2} and taking an approximation of the fluid volume fraction in terms of fluid pressure and volume change we obtain \eqref{eq3:mass}, below.

For transforming the momentum balance \eqref{p:momentum} we use the divergence theorem, the change of coordinates, we recall the relation between Cauchy and first Piola--Kirchhoff stress $\bP = J\bsigma\bF^\intercal$, and use again the divergence theorem to get \eqref{eq3:momentum}, below. The hyperelastic component to the stress is  
\begin{linenomath*}
\begin{equation}\label{eq:P}
\bP_e = \frac{\partial \Psi}{\partial \bF}|_{\widetilde{\Omega}},\end{equation}
\end{linenomath*}
and for a neo--Hookean material, equations \eqref{eq:FeFg}-\eqref{eq:P} imply that the effective stress tensor is 
\begin{linenomath*}
\begin{equation}\label{eq:Pactive}
\bP_e = J  (\mu \bB_e - \psi \bI) \bF^{-\intercal},
\end{equation}
\end{linenomath*}
where $\mu$ is the shear modulus, $\psi $ denotes a Lagrange multiplier,
and $\bB_e = \bF_e\bF_e^\intercal$ is the right Cauchy--Green deformation tensor associated with the intermediate
configuration    $\widetilde{\Omega}$.
Note that in this context, rather than material incompressibility of the mixture, the Lagrange multiplier $\psi$ enforces  a mass conservation kinematic constraint
of the form $J_e=1$, confirming that for a hyperelastically incompressible material 
all volume changes are solely due to growth factors
($J_g =  \det \bF_g$ is in turn irreversible) \cite{benamar05}. If $J_g>1$ the body undergoes growth. 

In summary, the equations of nonlinear poroelasticity with growth
(including momentum conservation of the solid-fluid mixture, fluid mass balance, and
incompressibility of elastic deformations)
read as follows
\begin{linenomath*}
\begin{subequations}
\begin{align}\label{eq3:momentum}
\rho\ddot\bu-\bDiv (\bP_e + \bP_f) & = \cero & \text{in } \Omega_0\times(0,t_{\mathrm{final}}], \\
\dot{\overline{\bigl(C_0 p+\alpha_{BW}[\tr\bC-d]\bigr)}} -\frac{1}{J}
\vDiv\biggl(\frac{\kappa}{\eta} J \bC^{-1} \vGrad p\biggr)  &= \ell & \text{in } \Omega_0\times(0,t_{\mathrm{final}}], \label{eq3:mass}\\
J_e & =1 & \text{in } \Omega_0\times(0,t_{\mathrm{final}}],\label{eq3:incompress}
\end{align}\end{subequations}\end{linenomath*}
where the stress exerted by the fluid (when pulled back to the reference configuration)
is simply $\bP_f = - \alpha_{BW}J  p  \bF^{- \intercal}$, and the divergence and gradient operators are now
understood with respect to the undeformed coordinates. In contrast with the model from
Section~\ref{sec:model}, now the fluid has a source that we relate to the growth process
and make it precise below. The divergence of displacement
accompanying the Biot--Willis constant has now been replaced by $\tr\bC-d$, where  $d$ is the spatial dimension.

It remains to characterise the growth deformation gradient. We suppose that growth
occurs due to a set of internal variables $\gamma_1,\gamma_2,\gamma_3$ acting on
local orthonormal directions $\bk_1$, $\bk_2$, $\bk_3$ in the undeformed configuration (and in 2D we consider
only $\gamma_1,\gamma_2$ acting on  $\bk_1$, $\bk_2$, respectively), and a general form for
the growth deformation tensor is simply
\begin{linenomath*}
\[\bF_g = \bI + \sqrt{\gamma_{1}}\, \bk_1\otimes\bk_1+ \sqrt{\gamma_{2}}\, \bk_2\otimes\bk_2+ \gamma_3 \bk_3\otimes\bk_3.
\]\end{linenomath*}
Assuming that a transversely isotropic growth due to a constant rate occurs on the plane $(\bk_1,\bk_2)$ and the
growth due to morphogen concentrations acts mainly on the direction $\bk_3$, we employ
the following specification for the growth factors $\gamma_i$:
\begin{linenomath*}
\begin{equation}\label{eq:gamma}
\gamma_{1} = \delta_1 t, \quad \gamma_{2} = \delta_1 t, \quad \gamma_3 = \delta_2t + \delta_3 \frac{m}{1+m^2},
\end{equation}\end{linenomath*}
where $\delta_1,\delta_2,\delta_3$ are positive constants and we recall that $m$ is the concentration of
mesenchymal cells in the form used in the active stress from \eqref{eq:active}.   In the 2D case, \eqref{eq:gamma}
will be replaced by
\begin{linenomath*}
\begin{equation}\label{eq:gamma2D}
\gamma_{1} = \gamma_{2} = \delta_2t + \delta_3 \frac{m}{1+m^2}.\end{equation}
\end{linenomath*}
As in \cite{giverso15,mascheroni18} a more general description that defines a functional relation between the rates of the growth factors and the exchange of mass implies a dependence on the solid volume fraction 
\begin{linenomath*}
\[ \sum_{i=1}^d\frac{\dot{\gamma_i}}{\gamma_i} = \frac{r_{s}}{\phi_s}.\]
\end{linenomath*}

Following \cite{hine16}, we can write a constitutive equation for the mass source
in \eqref{eq3:mass} depending on the constants $\gamma_i$ as follows
\begin{linenomath*}
\begin{equation}\label{eq:fluid-mass}
\ell = \ell_0( 2\delta_1+\delta_2 + \delta_3 \partial_m(\gamma_3) \partial_t m),
\end{equation}\end{linenomath*}
where $\ell_0$ is a constant and
where we stress that $\rho$ in this section refers to the density of the solid in the reference configuration.

The set of equations is complemented with suitable boundary and initial conditions, formulated
 in the reference configuration. In contrast with \eqref{bc:Gamma}-\eqref{bc:Sigma}, here
 we impose Robin conditions for the poroelasticity on the whole undeformed boundary
\begin{linenomath*}\begin{subequations}\begin{align}
(\bP_e + \bP_f)\bN_0 + \xi_R J  \bF^{- \intercal} \bu & = \cero \qon (\Gamma_0\cup \Sigma_0)\times(0,t_{\text{final}}], \label{bc:robin}\\
\frac{\kappa}{\eta}J \bC^{-1}\vGrad p \cdot\bN_0 & = 0 \qon \Gamma_0\times(0,t_{\text{final}}], \qquad p = p_D \qon \Sigma_0\times(0,t_{\text{final}}],\label{bc:pressure}
\end{align}\end{subequations}\end{linenomath*}
where $\xi_R$ is the stiffness of the springs attached to $\partial\Omega_0$.

Also the chemotaxis system is recast in the reference domain. In particular, the advection terms are absorbed
by the material derivatives (since we are assuming that all species are transported by the whole mixture and can diffuse in both phases) and the diffusion coefficients are modified by the Piola transformation. As in
\eqref{eq3:mass}, the divergence of displacement in the sink term from \eqref{eq:f} now is replaced by $\tr\bC-d$. The transformation also involves Reynolds transport theorem, divergence theorem, and localisation. 
Assuming that
$D_m,D_f,D_b$ are isotropic, we have
\begin{linenomath*}
\begin{subequations}
 \begin{align}
  \dot{m} - 	\frac{1}{J}\vDiv\bigl( D_mJ \bC^{-1} \, \vGrad m - \alpha m \exp(-\gamma m) J \bC^{-1} \vGrad f \bigr)
&= 0 & \text{in } \Omega_0\times(0,t_{\mathrm{final}}],	\label{eq2:m}\\
	\label{eq2:e}
\dot{e} -[ \kappa_1 w(\bX,t) h_1(m)+\kappa_2h_2(m)] (1-e)+[1-h_1(m)](\kappa_3+\kappa_4b)e & =0 & \text{in } \Omega_0\times(0,t_{\mathrm{final}}],\\
\label{eq2:f}
  \dot{ f}  -
	\frac{1}{J}\vDiv\bigl( D_fJ \bC^{-1}\, \vGrad f \bigr)
- \kappa_{F}e + \delta_{F}f + \xi_f (\tr\bC-d) & = 0 & \text{in } \Omega_0\times(0,t_{\mathrm{final}}],\\
\label{eq2:b}
  \dot{b}  -
	\frac{1}{J}\vDiv\bigl( D_bJ \bC^{-1}\, \vGrad b \bigr)
- \kappa_{B}h_3(m)m + \delta_{B}b& = 0 & \text{in } \Omega_0\times(0,t_{\mathrm{final}}],
\end{align}\end{subequations}\end{linenomath*}
where divergence and gradients are taken with respect to the undeformed configuration, and the zero-flux boundary conditions \eqref{bc:m} are changed accordingly
\begin{linenomath*}
\begin{align*}
\bigl(D_mJ \bC^{-1} \vGrad m - \alpha m \exp(-\gamma m)J \bC^{-1} \vGrad f \bigr) \cdot \bN_0 & = D_fJ \bC^{-1} \vGrad f \cdot\bN_0 \\
& =  D_mJ \bC^{-1} \vGrad m  \cdot\bN_0 \\&= 0 \qquad \text{on $\partial\Omega_0\times(0,t_{\text{final}}]$}.\end{align*}
\end{linenomath*}
\section{Numerical tests}\label{sec:results}

\begin{figure}[!t]
\begin{center}
\includegraphics[width=0.245\textwidth]{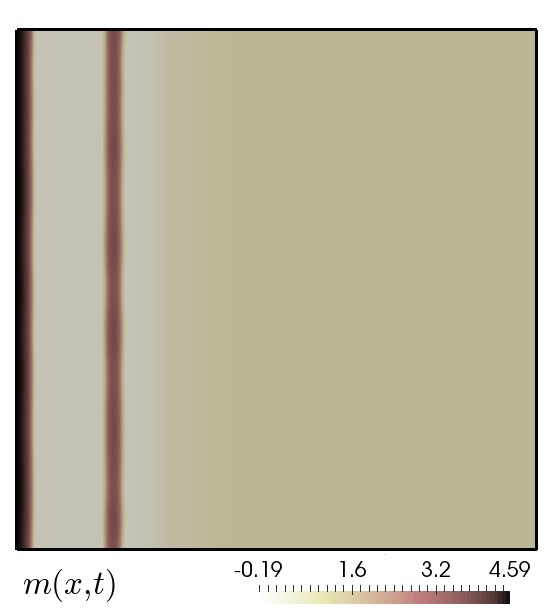}
\includegraphics[width=0.245\textwidth]{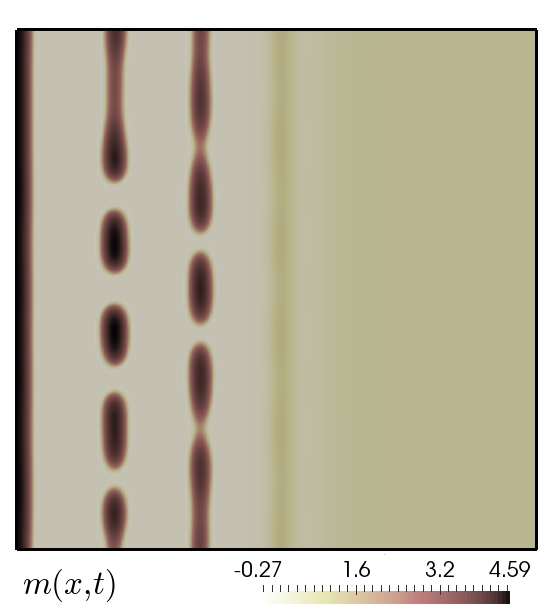}
\includegraphics[width=0.245\textwidth]{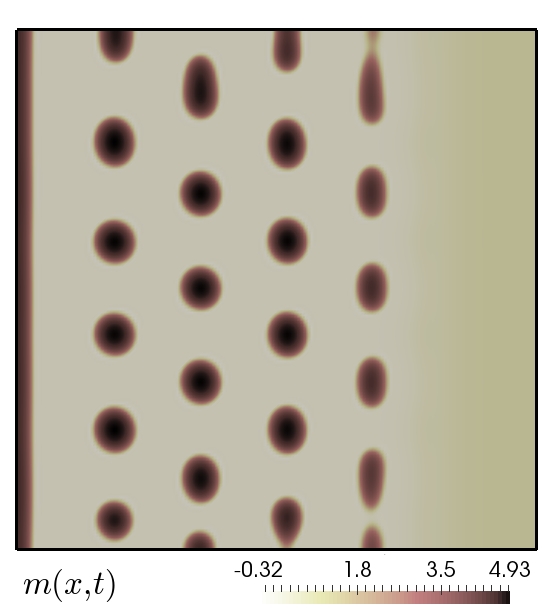}
\includegraphics[width=0.245\textwidth]{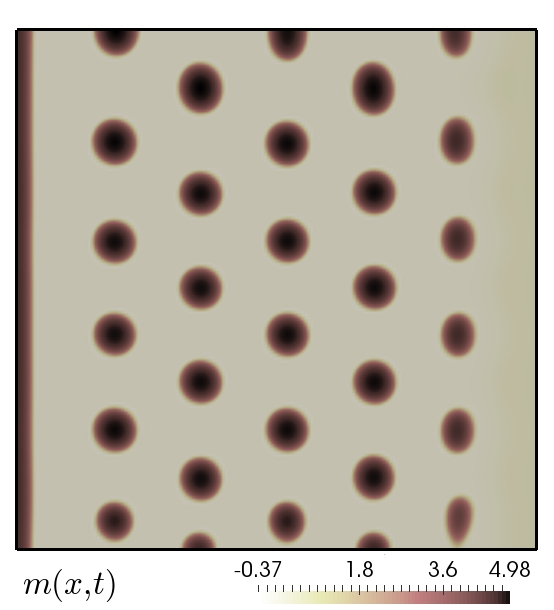}\\
\includegraphics[width=0.245\textwidth]{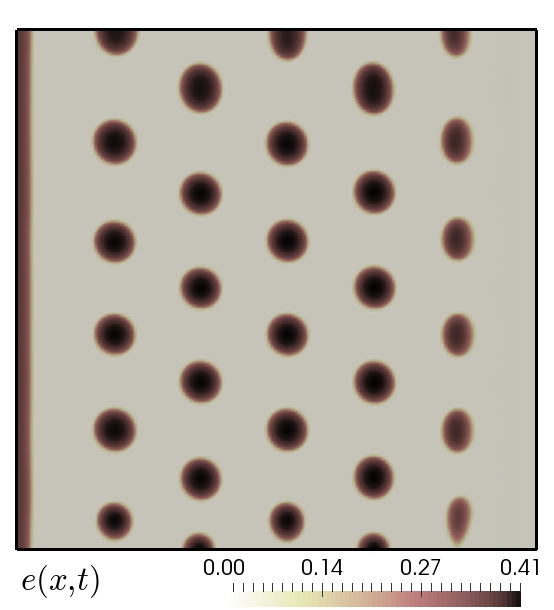}
\includegraphics[width=0.245\textwidth]{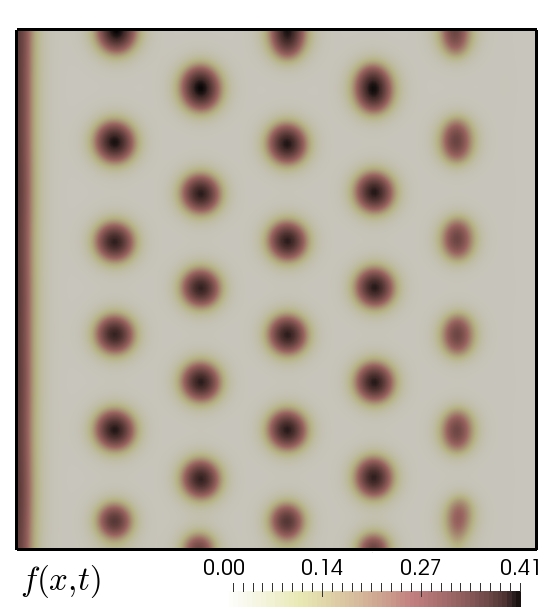}
\includegraphics[width=0.245\textwidth]{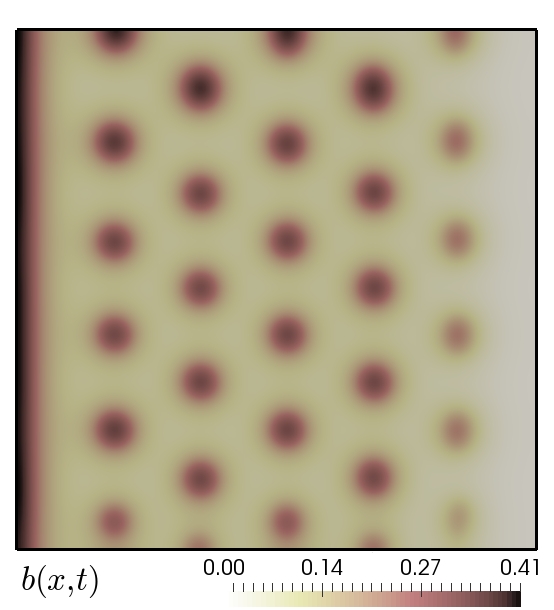}
\end{center}
\caption{Evolution of the mesenchymal cell concentration at $t=160,280,400,520$, under no deformation (top panels); and epithelium, FGF, and BMP concentrations at $t=520$ (bottom).}
\label{fig:ex01a}
\end{figure}

\subsection{Discretisation and implementation.}
	We have employed mixed finite element methods for the spatial discretisation of systems \eqref{eq:m}-\eqref{eq:active} and \eqref{eq3:momentum}-\eqref{eq2:b}. We use quasi-uniform triangular/tetrahedral meshes in all cases. For the first system described in Section~\ref{sec:model}, since the linear poroelasticity equations are away from the nearly incompressibility limit (here the Poisson ratio is $\nu = 0.4$), it suffices with using piecewise linear and continuous approximation of solid displacements (as well for all other unknowns). For the system proposed in Section~\ref{sec:growth}, we employ the so-called Taylor-Hood mixed approximation of solid displacements and solid pressure composed by piecewise quadratic and continuous functions for each displacement component combined with piecewise linear and overall continuous functions for the solid pressure (this pair is known to satisfy an adequate inf-sup condition for the linearised hyperelasticity equations, see, e.g., \cite{braess05}); and the fluid pressure as well as the chemical species concentrations are approximated with piecewise linear and continuous polynomials. For all problems, an outer fixed-point algorithm decouples the chemotaxis system from the poromechanics. An implicit-explicit method is employed to generate a full discretisation of the chemotaxis system and of the linear poroelasticity equations from Section~\ref{sec:model}, whereas for the poroelastic-growth system \eqref{eq3:momentum}-\eqref{eq3:incompress} with boundary conditions \eqref{bc:robin}-\eqref{bc:pressure}, an inner Newton--Raphson method is used to solve the corresponding nonlinear equations written in implicit form. All routines for the 2D tests have been implemented in the open-source finite element library FEniCS \cite{alnaes15}. The solution of all nonlinear systems was carried out with Newton's method using a tolerance of 1E-7 on either the  relative or absolute $\ell^2-$norm of the vector residuals, and the linear systems on each step were solved with the distributed direct solver MUMPS. For the 3D case we used the Firedrake library \cite{firedrake} due to its facility in handling block preconditioners. We provide details regarding the preconditioners in Section \ref{sec:preconditioners}. The nonlinear tolerances are taken identical to the 2D case. For the linear solvers, we use a GMRES with the nested Schur preconditioner recently proposed in \cite{barnafi22}, with an absolute tolerance of 1E{-12} and a relative one of 1E{-2}. The use of such a big relative tolerance yields a simplified inexact-Newton method \cite{inexnewt} that results in more nonlinear iterations that require much less time.

	\subsection{Efficient preconditioners.}\label{sec:preconditioners}
Given that we are using a Newton solver, an efficient preconditioner is fundamental to obtain computational times that are feasible for the 3D case. The splitting strategy we are considering allows us to consider the poroelastic and chemotaxis problems separately, where each of them consideres a different preconditioning strategy. In what follows we give the details of the preconditioners used for the solution of the tangent operator for each physics.

\begin{itemize}
	\item {\bf Chemotaxis.} This problem considers four similar building block physics, consisting essentially of three parabolic problems $(m,f,b)$ and one algebraic constraint $(e)$. We consider an additive block solver, meaning that we use a block-wise Jacobi preconditioner with only the diagonal block of each variable. At the block level, we consider the action of an AMG preconditioner for the parabolic problems and the action of a Jacobi preconditioner for the algebraic one. 
	\item {\bf Poromechanics.} The poromechanics block is more difficult, which is reflected in the complexity of the preconditioner under consideration. It is based on the block preconditioner proposed in \cite{barnafi22} for large-strain  poromechanics, where we use a lower Schur complement block factorisation with the fields $\bu$ and $(\psi, p)$. For such a  block we consider the action of an AMG preconditioner, whereas for the corresponding Schur complement block we consider instead a sparse representation given by an ad-hoc extension of the fixed-stress splitting scheme proposed in \cite{borregales2018}. For this, we add two stabilisation terms given by 
	\begin{linenomath*}
		\[ \int_{\Omega_0} \beta_s\psi \psi^* + \beta_f pp^*\dx, \]
		\end{linenomath*}
where we have denoted with $(\cdot)^*$ the test function corresponding to each variable and  have used the values $\beta_s = 1.0$ and $\beta_f=0.1$ \cite{borregales2018}. As this approach yields a sparse operator that approximates the Schur complement, we simply use the action of an AMG preconditioner based on it.
\end{itemize}

We highlight that the proposed preconditioner represents an improvement over \cite{barnafi22}, as we employ only preconditioner actions on each block, which allows us to use a GMRES linear solver instead of a flexible GMRES (which can be substantially more demanding in terms of computational cost).

\begin{figure}[!t]
\begin{center}
\includegraphics[width=0.245\textwidth]{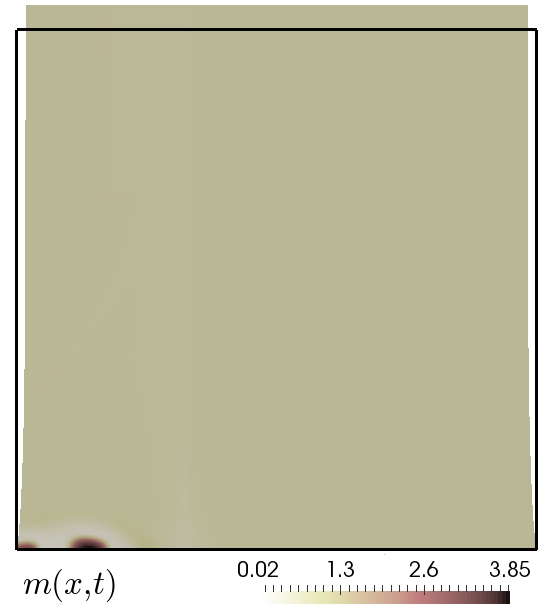}
\includegraphics[width=0.245\textwidth]{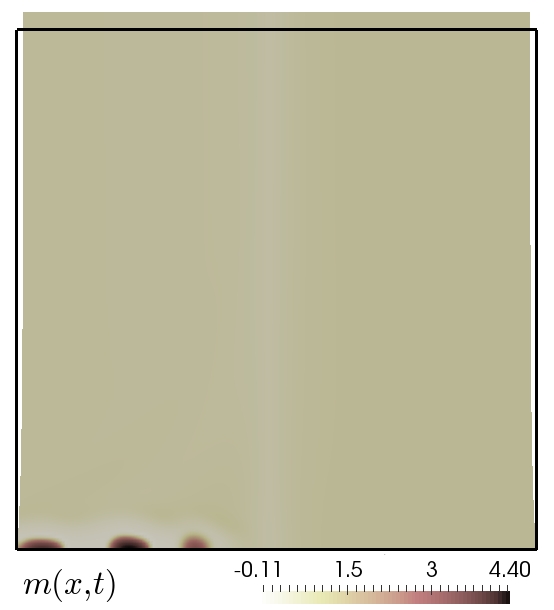}
\includegraphics[width=0.245\textwidth]{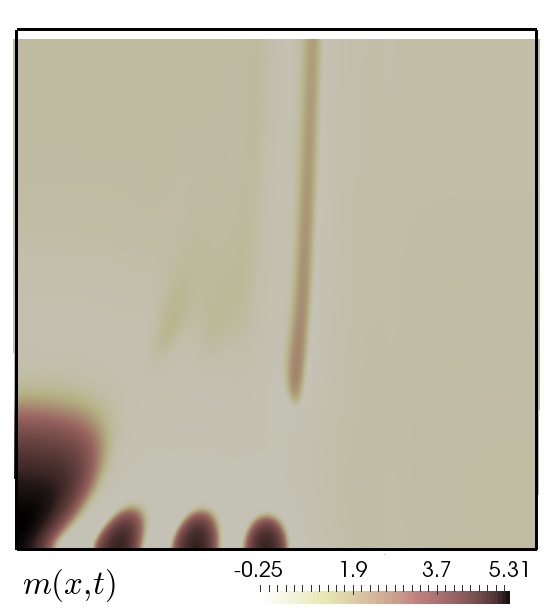}
\includegraphics[width=0.245\textwidth]{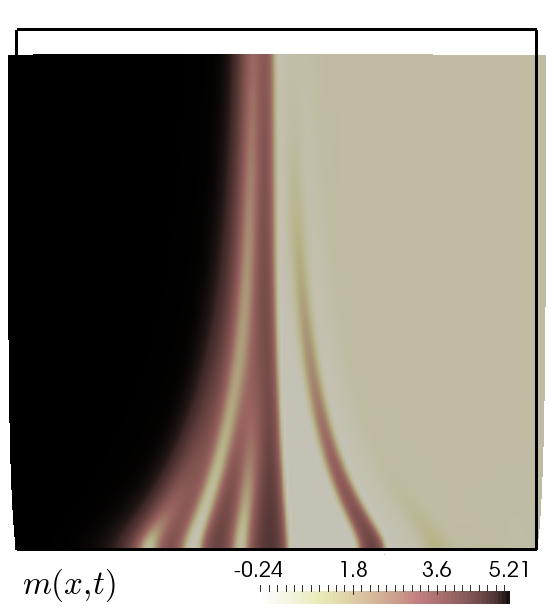}\\
\includegraphics[width=0.245\textwidth]{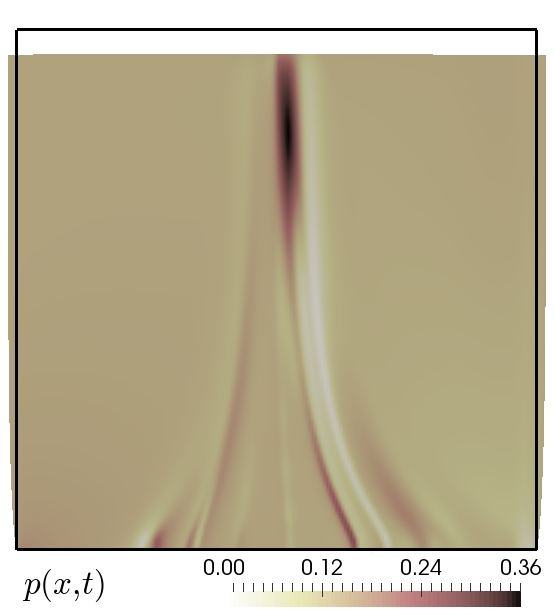}
\includegraphics[width=0.245\textwidth]{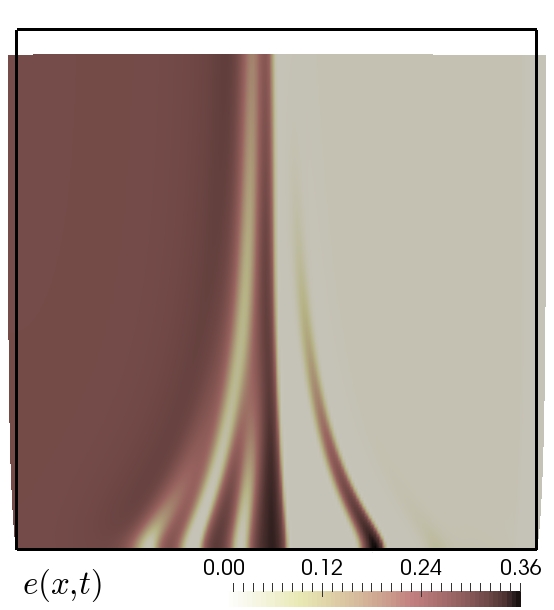}
\includegraphics[width=0.245\textwidth]{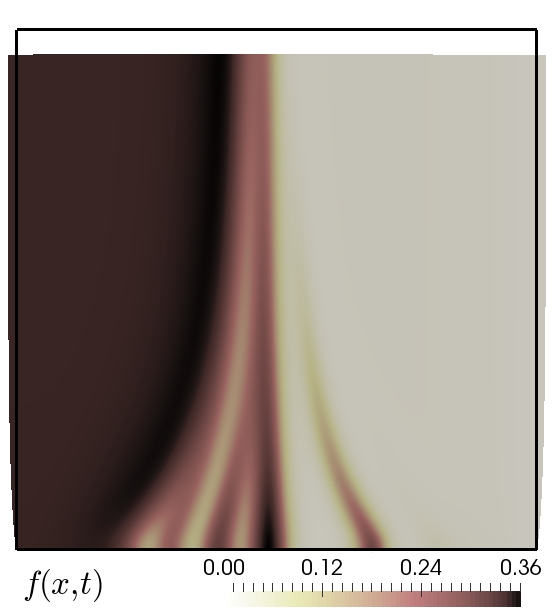}
\includegraphics[width=0.245\textwidth]{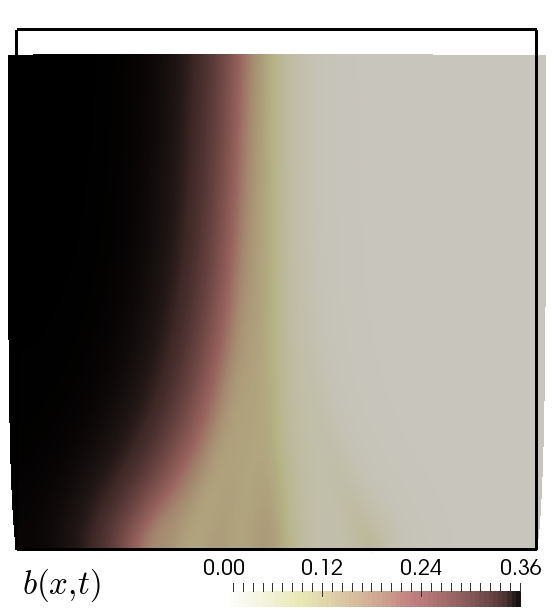}
\end{center}
\caption{Evolution of the mesenchymal cell concentration at $t=160,280,400,520$ under periodic traction applied on the top edge of the domain (top); and snapshots of fluid pressure, epithelium, FGF, and BMP at $t=520$ (bottom row).}
\label{fig:ex01b}
\end{figure}

\subsection{Suppressed solid motion vs. periodic boundary traction.}
In order to emphasise the effect of the stretch and deformation on the formation of patterns we compare the concentration of mesenchymal cells in two cases, one where the motion of the domain is suppressed (by setting $\tau = 0$), and another scenario where the bottom of the domain, $(\Gamma)$, is kept clamped, the top of the domain undergoes a (relatively small) periodic traction, and the vertical faces are kept with zero traction. In this case we employ the linear poroelastic model from Section~\ref{sec:model}, and a zero fluid pressure flux is imposed on $\Sigma=\partial\Omega\setminus\Gamma$. The spatial domain is a square $\Omega = (0,20)^2$ discretised into 58482 triangular elements, and the model parameters that are modified with respect to those in Table~\ref{tab:LAparams} adopt the values
\begin{linenomath*}\[m_0 =0.75 ,\quad  D_m = 0.03, \quad \tau = \xi_f = 0.001.\]\end{linenomath*}
The traction is $\bt = (0,s_0\sin(\pi t/\hat{t}))^\intercal$, with magnitude $s_0 = 1500$ and period $\hat{t} = 320$. We employ a fixed time step $\Delta t=0.2$ and advance the system until $t=520$. The outcome of these tests is collected in Figure~\ref{fig:ex01a} and Figure~\ref{fig:ex01b}, where we observed that slight modifications in the poromechanics result in quite drastic alterations of the clustering of mesenchymal cells.  For instance, in Figure~\ref{fig:ex01b} we see that on the bottom edge of the domain (where displacements are zero during the whole simulation), the relatively large stretches applied elsewhere on the body do not break the formation of dotted-shaped patterns at early stages, but they are swept away as time advances.

 \begin{figure}[!t]
\begin{center}
\includegraphics[width=0.245\textwidth]{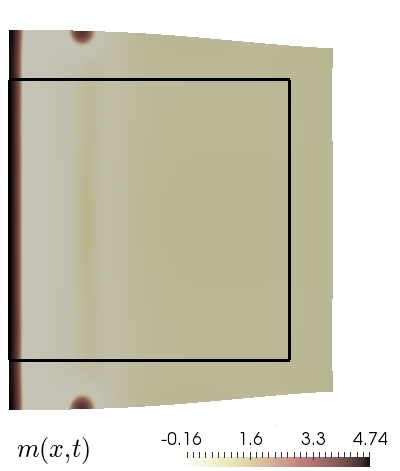}
\includegraphics[width=0.245\textwidth]{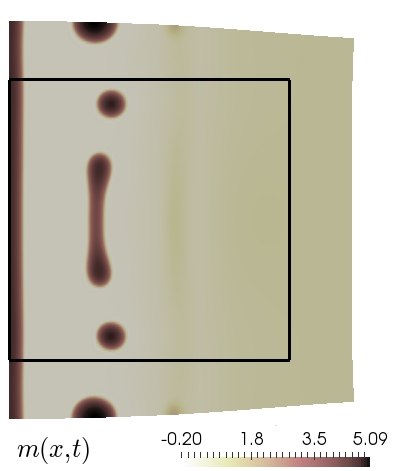}
\includegraphics[width=0.245\textwidth]{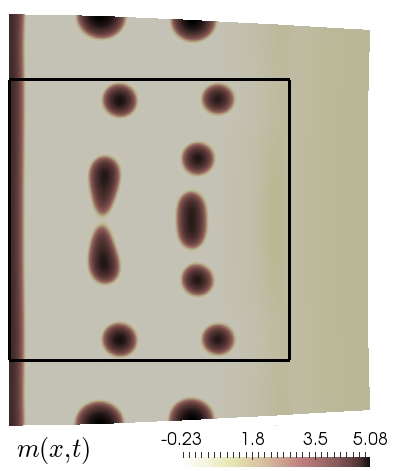}
\includegraphics[width=0.245\textwidth]{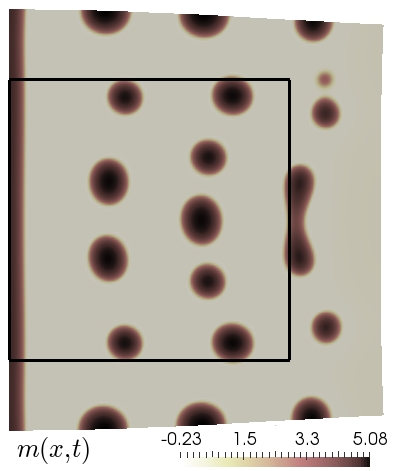}\\
\includegraphics[width=0.195\textwidth]{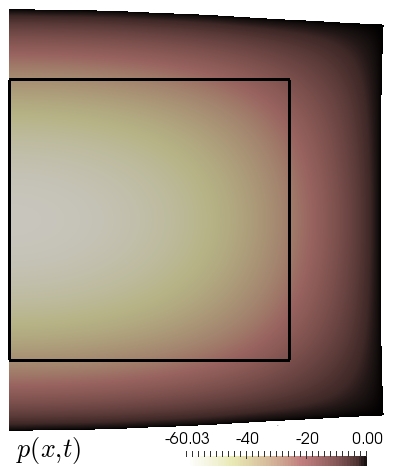}
\includegraphics[width=0.195\textwidth]{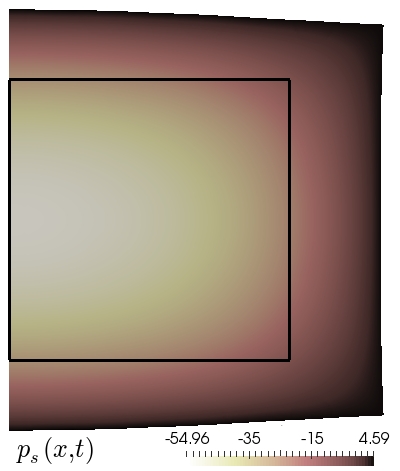}
\includegraphics[width=0.195\textwidth]{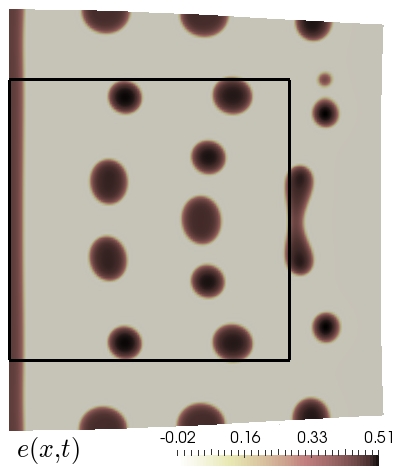}
\includegraphics[width=0.195\textwidth]{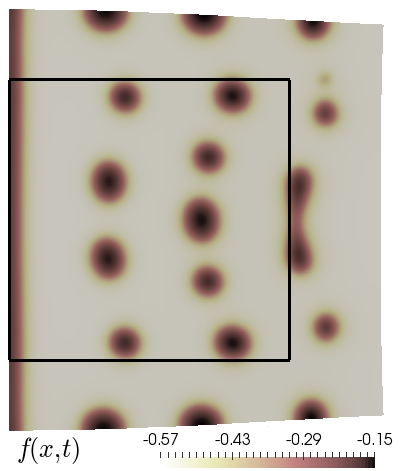}
\includegraphics[width=0.195\textwidth]{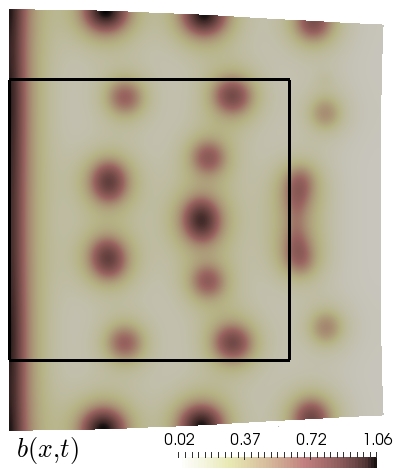}
\end{center}
\caption{Evolution of the mesenchymal cell concentration under finite growth using the formulation \eqref{eq:gamma2D}, snapshots are taken at $t=160,280,400,520$ (first row). Bottom: plots at $t=520$ of fluid pressure, solid pressure, epithelium, FGF, and BMP. Mechano-chemical coupling governed by $\xi_f=0.5$.}
\label{fig:ex02a}
\end{figure}

\subsection{Finite growth in 2D.} Next we maintain the domain and most model parameters as in the previous case and focus on the coupled model defined in Section~\ref{sec:growth}. Determining the total first Piola-Kirchhoff stress from \eqref{eq:Pactive} and the growth specification \eqref{eq:gamma2D} (with $\delta_2 =0.15$, $\delta_3= 0.045$), we modify also the boundary treatment and consider that $\Gamma_0$ is the left edge of the domain whereas the remainder of $\partial\Omega_0$ is $\Sigma_0$. On $\Gamma_0$ we impose zero normal displacement $\bu\cdot\bN_0 = 0$ and a Robin boundary condition with constant $\xi_R = 0.0001$ is set on $\Sigma_0$. The shear modulus is $\mu = 4$, the solid permeability is augmented to $\kappa = 0.001$, and we take $\xi_f=0.3$. The fluid mass source is specified by \eqref{eq:fluid-mass} with $\ell_0 = 0.001$.

\begin{figure}[!t]
\begin{center}
\includegraphics[width=0.245\textwidth]{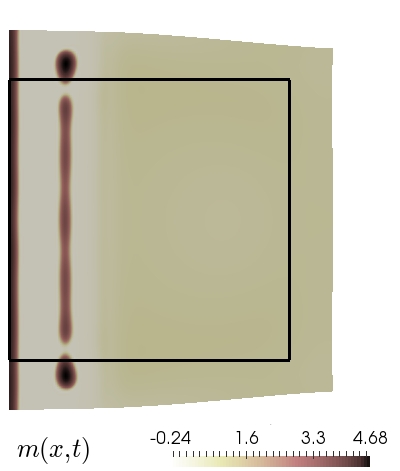}
\includegraphics[width=0.245\textwidth]{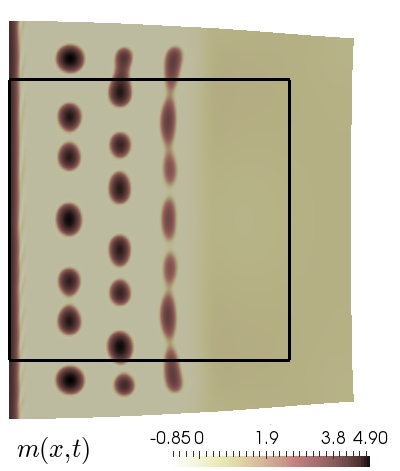}
\includegraphics[width=0.245\textwidth]{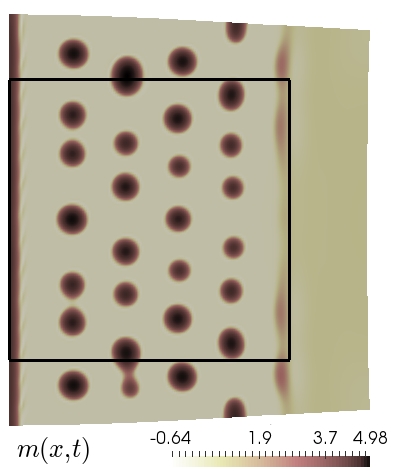}
\includegraphics[width=0.245\textwidth]{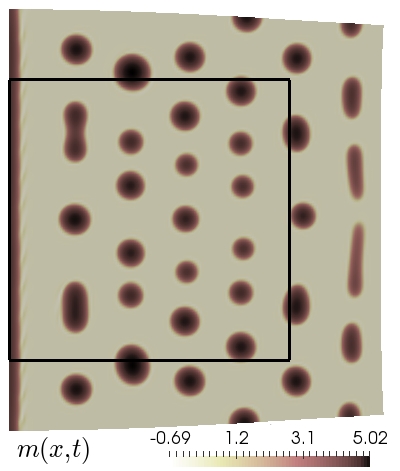}\\
\includegraphics[width=0.195\textwidth]{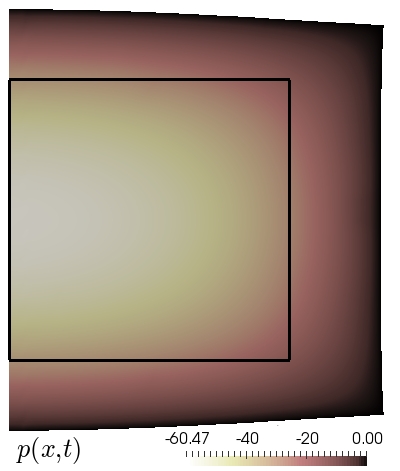}
\includegraphics[width=0.195\textwidth]{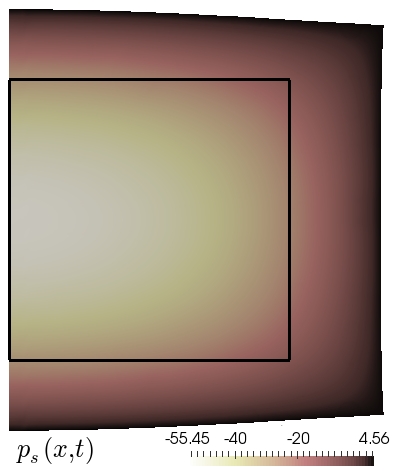}
\includegraphics[width=0.195\textwidth]{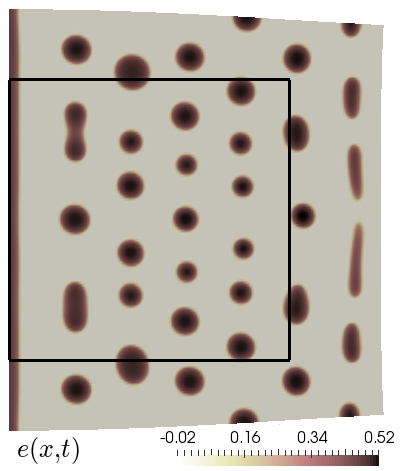}
\includegraphics[width=0.195\textwidth]{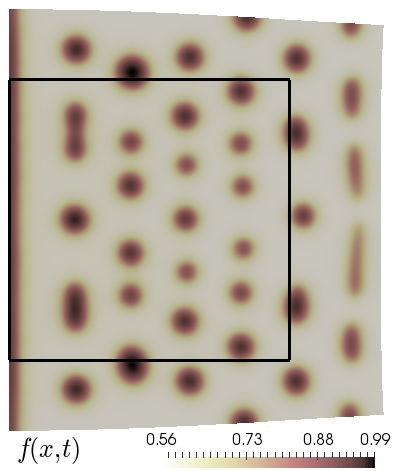}
\includegraphics[width=0.195\textwidth]{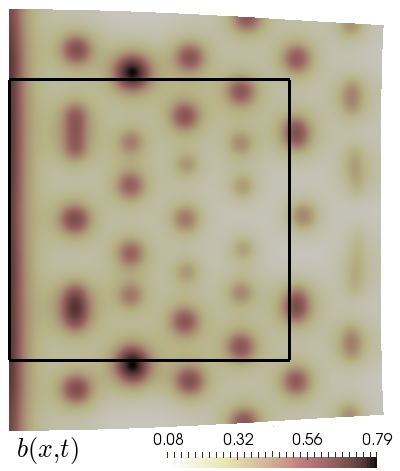}
\end{center}
\caption{Evolution of the mesenchymal cell concentration under
finite growth using the formulation \eqref{eq:gamma2D}, snapshots are taken at $t=160,280,400,520$ (top). Second row: plots at $t=520$ of fluid pressure, solid pressure, epithelium, FGF, and BMP. Simulations using $\xi_f=-0.5$. }
\label{fig:ex02b}
\end{figure}

The panels in Figure \ref{fig:ex02a} display transients of the growth of the initial domain and we plot the concentration of mesenchymal cells on the deformed configuration. We can observe that the patterns follow a different alignment than those formed in Figure~\ref{fig:ex01a}, in particular disrupting the symmetry and the size distribution of the dotted-shaped patterns.  The plots in the bottom row show the fluid and solid pressure together with the remaining chemical concentrations at the final time. Note that both solid and fluid pressures accumulate near $\Sigma_0$. In addition, we note that the values of FGF at later times are negative since the last term on the left-hand side of \eqref{eq2:f} represents a sink of $f$ as the domain only grows (and $\tr\bC-d$ is non-negative) throughout the simulation. We then change the sign of $\xi_f=-0.3$ and obtain positive FGF concentrations as well as similar patterns in all fields (but exhibiting a higher wave-number than in the previous case). The results from these tests are shown in Figure \ref{fig:ex02b}.

\subsection{Finite growth in 3D.} 
\begin{figure}[!t]
\begin{center}
\includegraphics[width=0.245\textwidth]{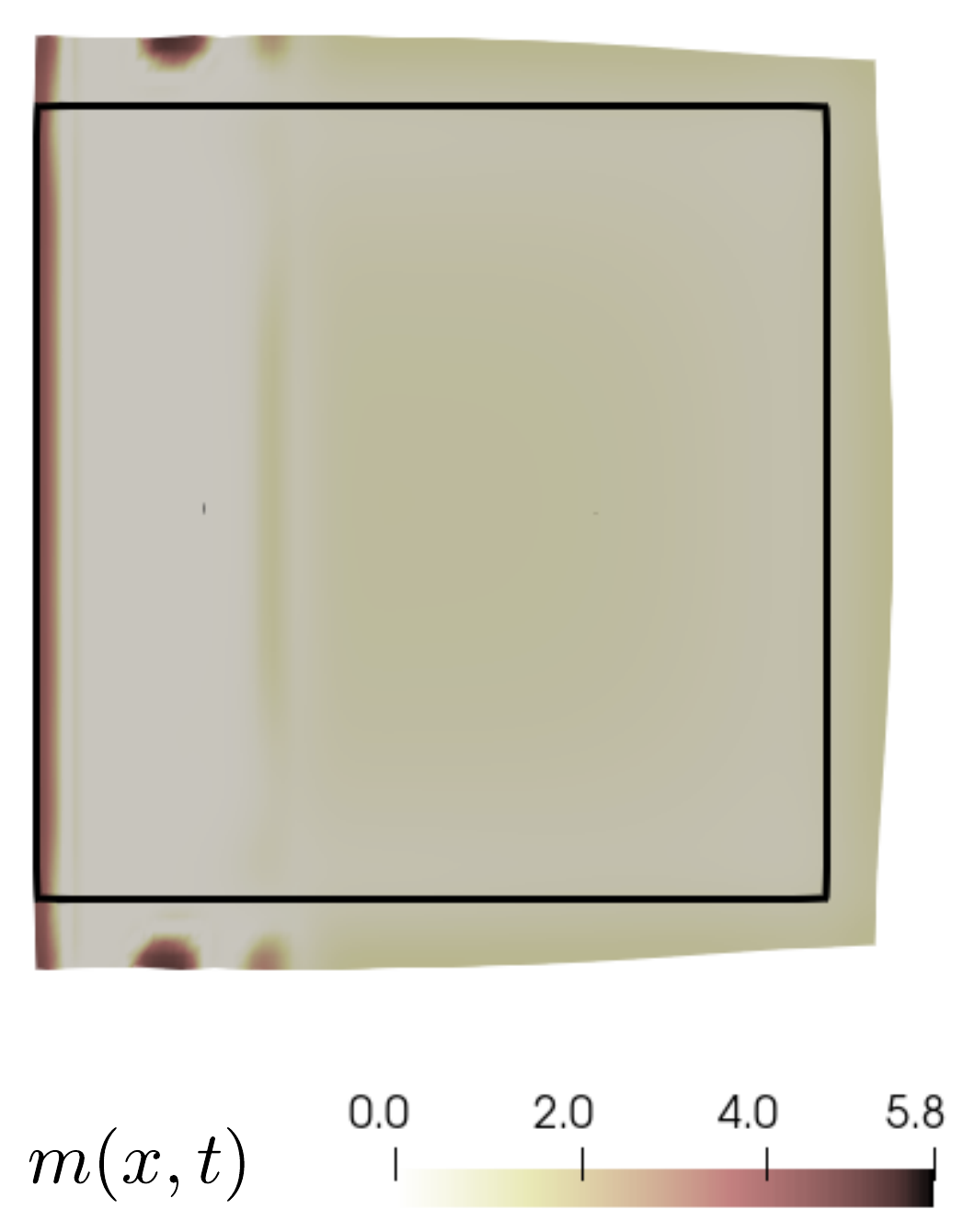}
\includegraphics[width=0.245\textwidth]{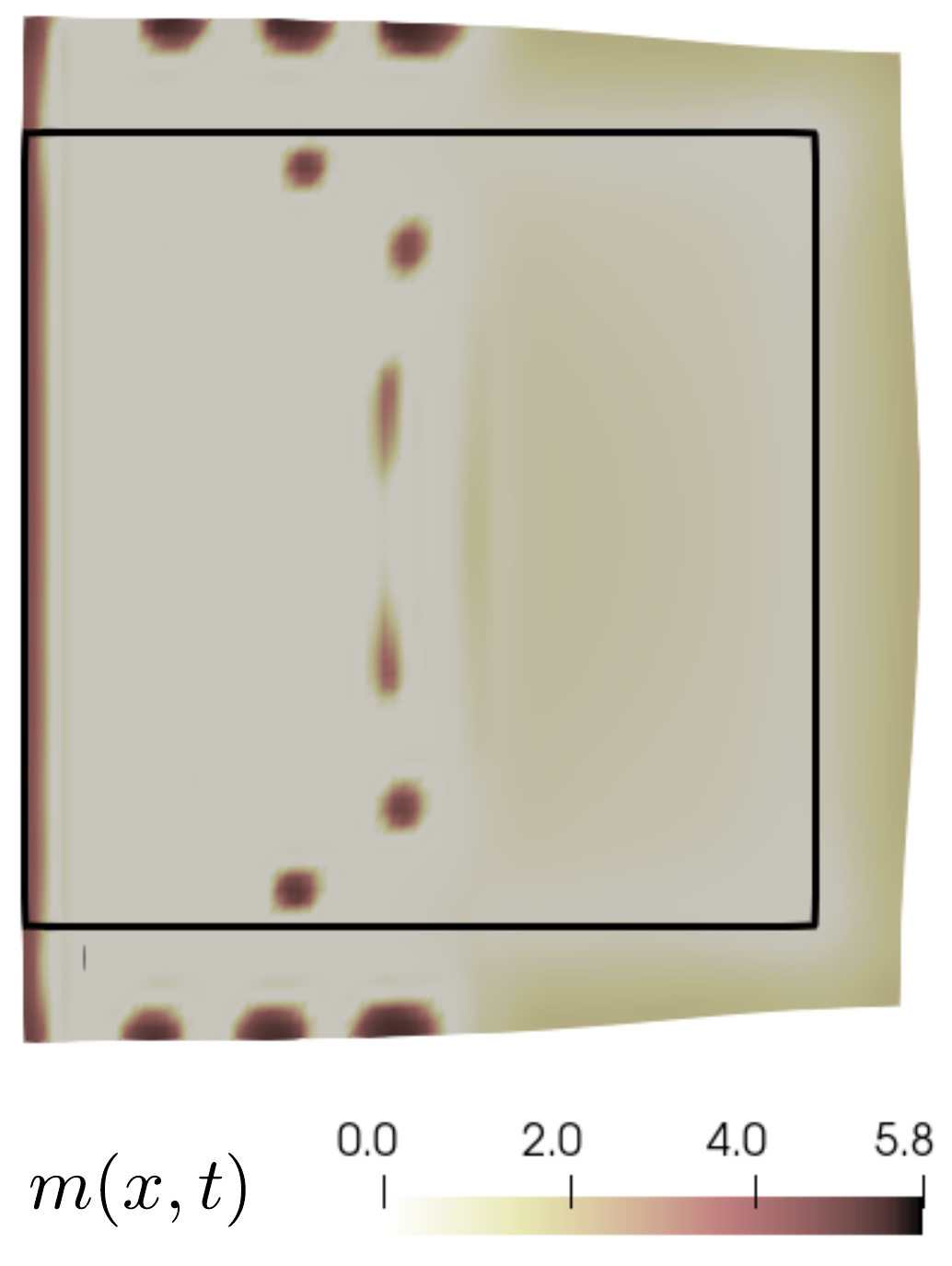}
\includegraphics[width=0.245\textwidth]{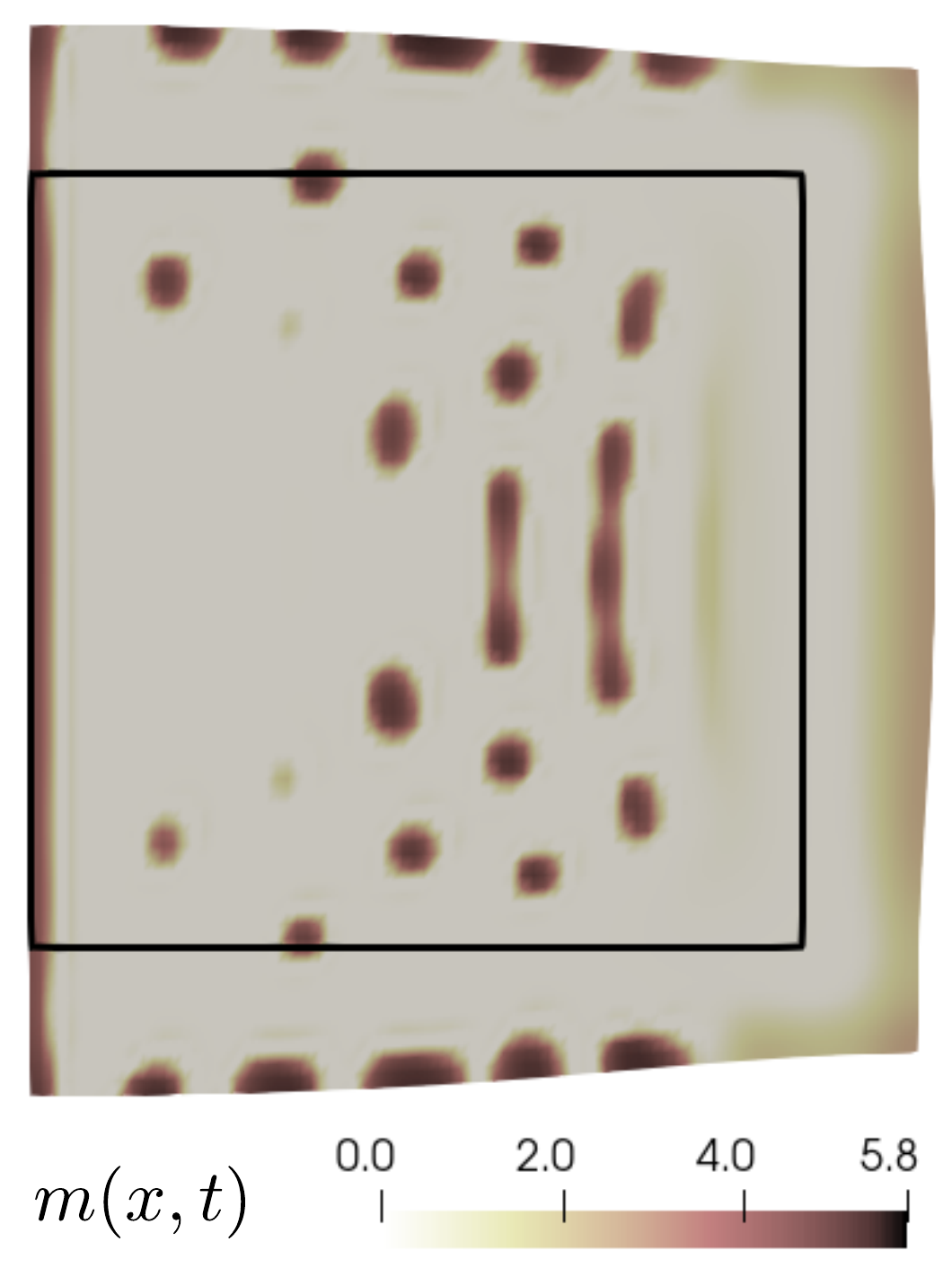}
\includegraphics[width=0.245\textwidth]{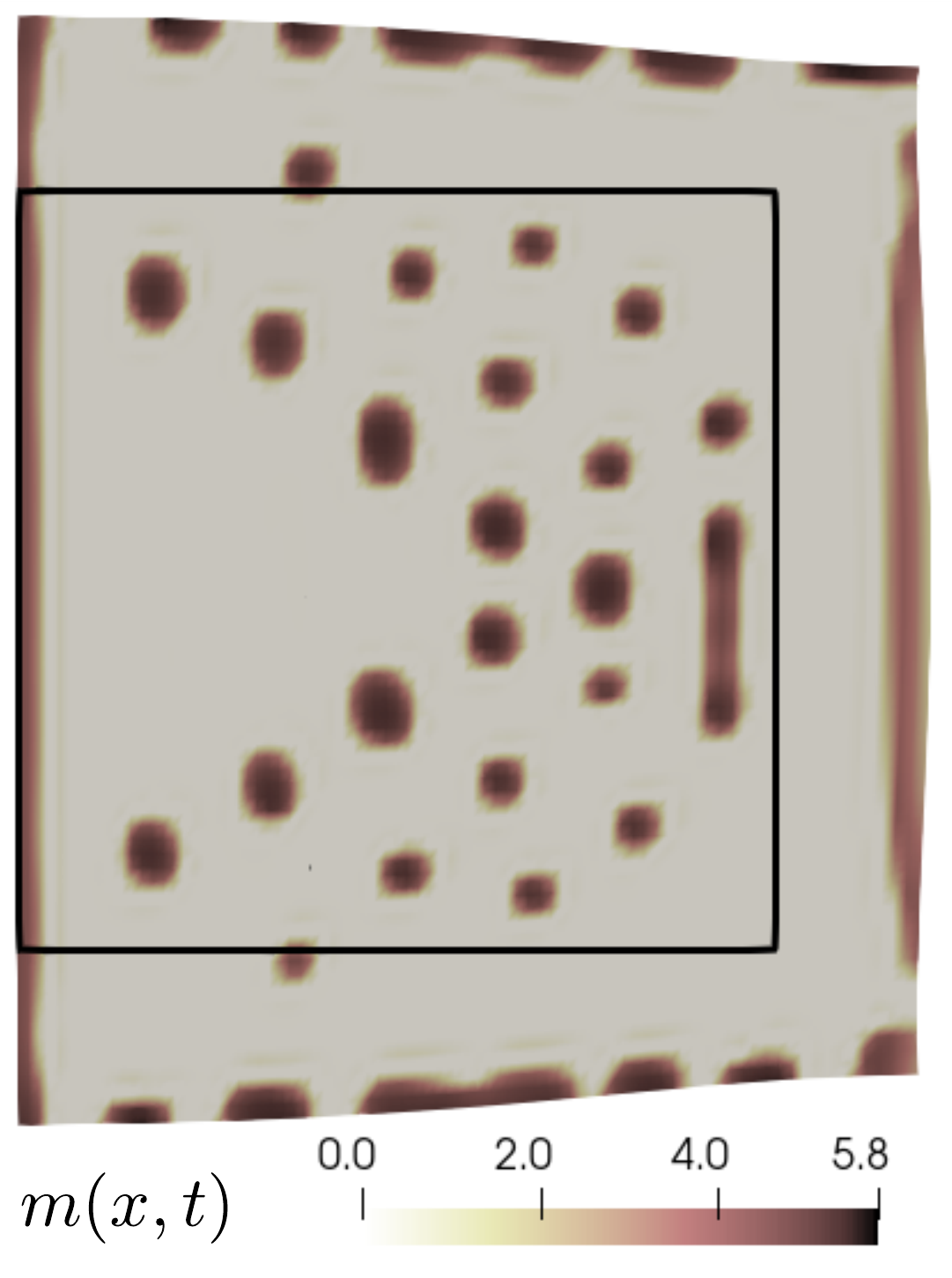}\\
\includegraphics[width=0.495\textwidth]{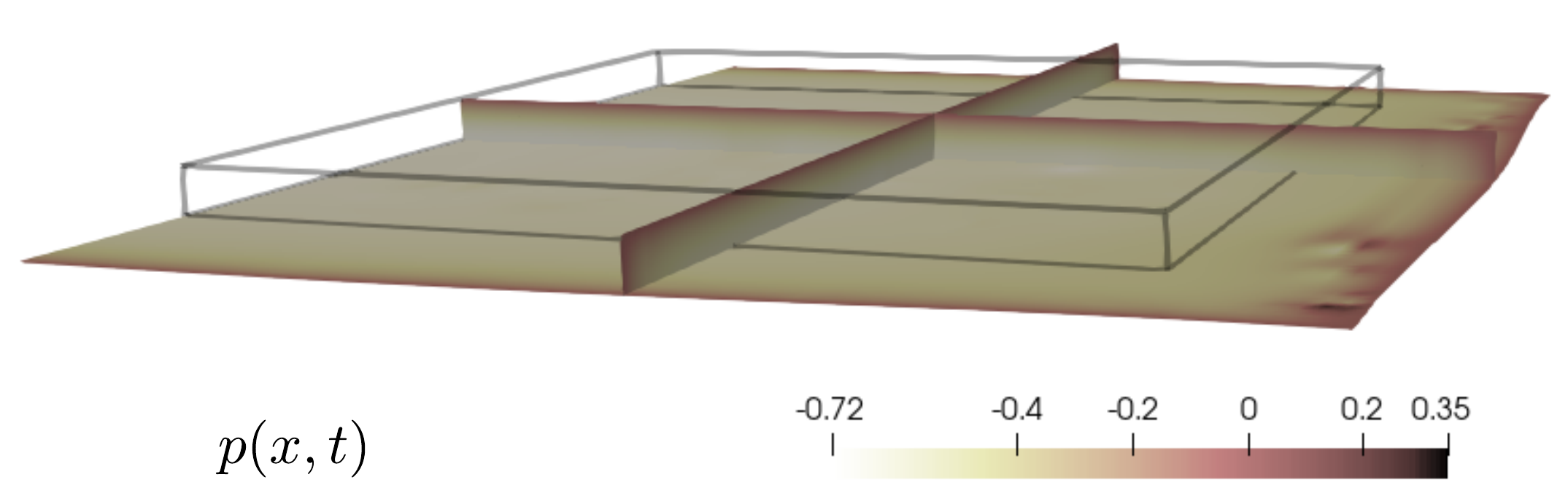}
\includegraphics[width=0.495\textwidth]{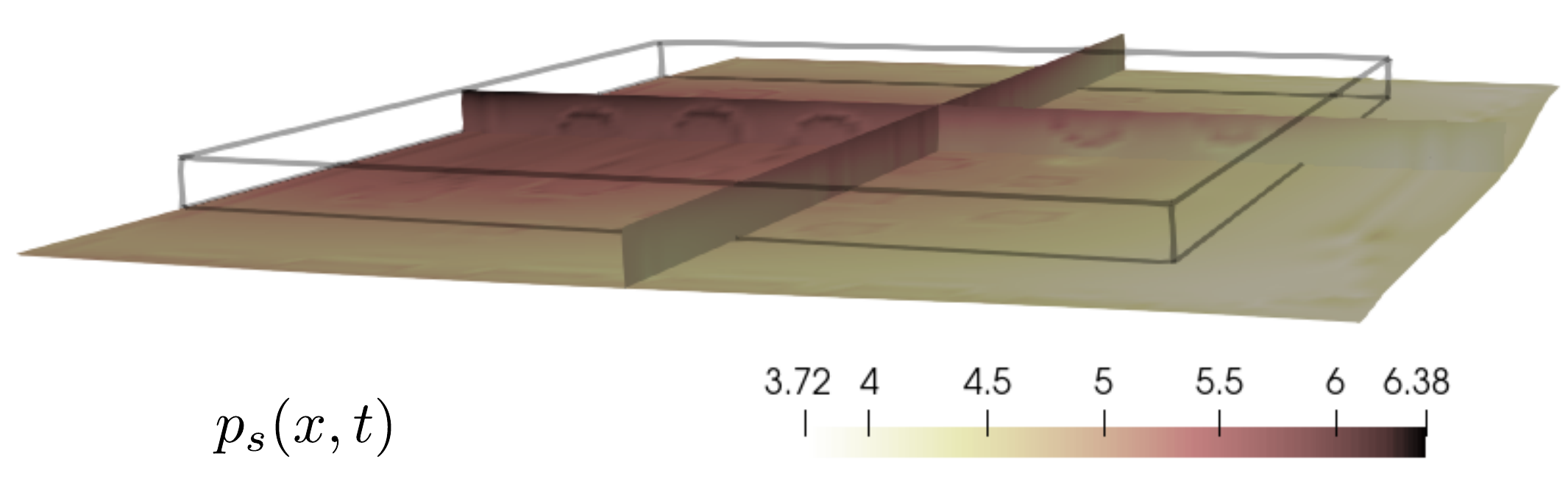}\\
\includegraphics[width=0.495\textwidth]{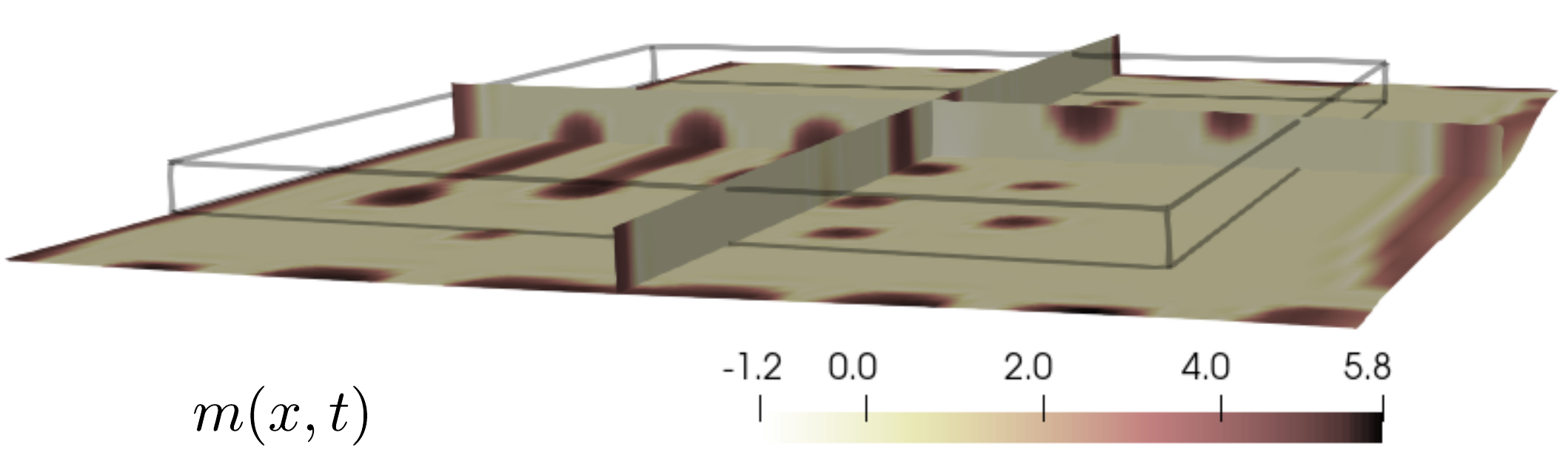}
\includegraphics[width=0.495\textwidth]{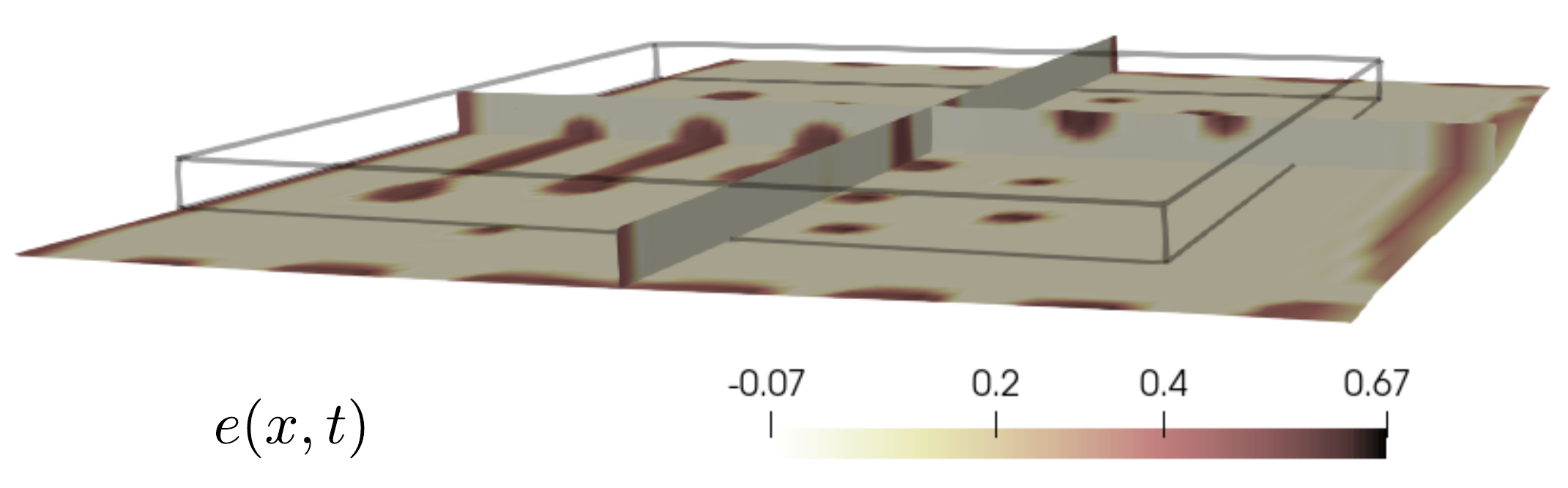}\\
\includegraphics[width=0.495\textwidth]{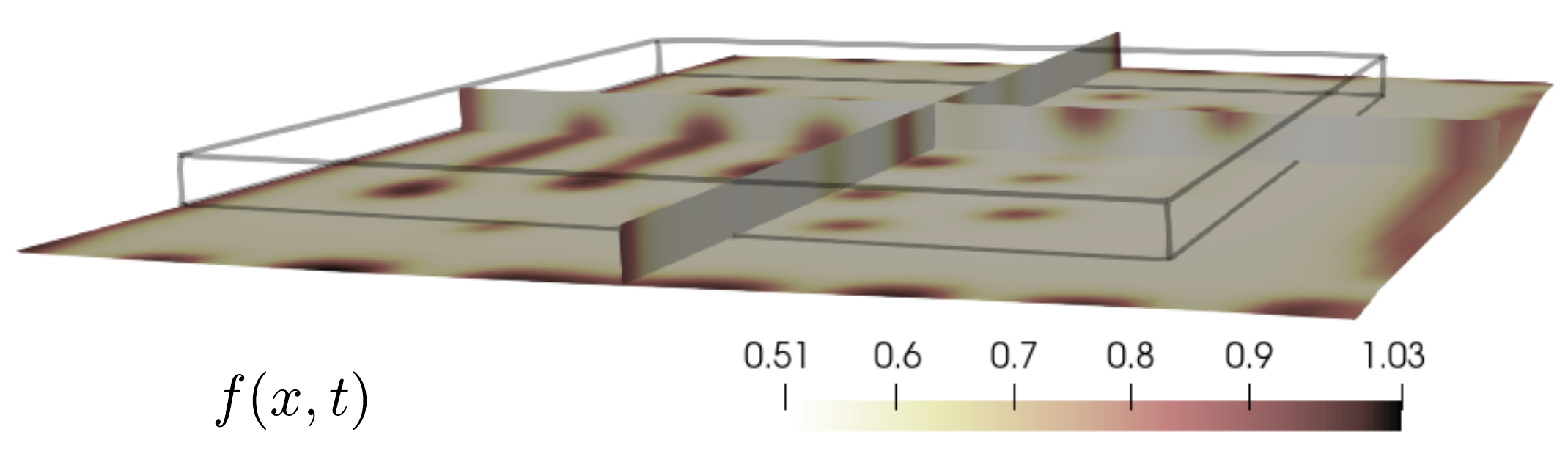}
\includegraphics[width=0.495\textwidth]{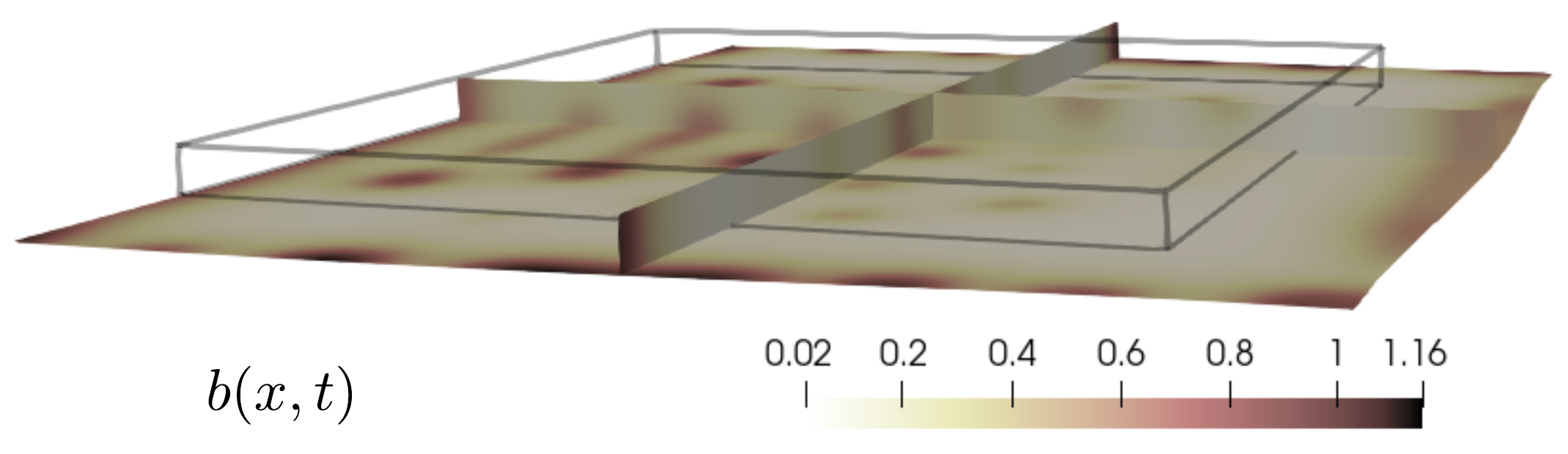}
\end{center}
\caption{Evolution of the mesenchymal cell concentration under finite growth using the formulation \eqref{eq:gamma}, snapshots are taken at $t=160,280,400,520$ (top). We also display the fluid and solid pressures (second row), the mesenchymal and epithelium concentrations (third row), and the FGF and BMP (fourth row).}
\label{fig:ex04}
\end{figure}

To close this section we proceed to extend the previous example to the 3D case. The undeformed domain is still a simple geometry $\Omega_0 = (0,20)\times (0,20)\times(0,1)$, where to obtain an accurate solution we used 140 elements per side in the $X_1$ and $X_2$ directions, and instead considered only 8 in the $X_3$ direction. This results in roughly 10 million degrees of freedom. The sides at $X_1=0$ and $X_3=0$ constitute $\Gamma_0$, where zero normal displacements are considered, and on the remainder of the boundary we set Robin conditions as above, but using $\xi_R=0.001$. The growth factors are considered as in \eqref{eq:gamma} with the specification $\delta_1=1/52$, $\delta_2=0$, $\delta_3=  1/104$. The results are shown in Figure~\ref{fig:ex04}, where we plot the evolution of the chemicals in the deformed configuration, and where the outline of the undeformed domain is also shown for visualisation purposes. We can see that as time progresses, the transversally isotropic growth depending only on time occurs mainly in the plane $X_1X_2$ as in the 2D case. 
Most interestingly, the mesenchymal concentration exhibits new complex structures. From the upper view, an arrow-like structure can be seen, where from higher concentrations appear more distant when going from right to left. From the slice view below, it can be seen how these distant concentrations are instead connected within the geometry, and we stress that this type of patterns can only be described by the 3D model. Still, these patterns are not representative of the actual primordial patterns, which are composed by independent buds, so this aspect still requires further investigation.

\section{Concluding remarks}\label{sec:concl}
We have proposed a model for describing the interaction between poroelastic deformations and chemotactical displacement of mesenchymal cells, FGF, BMP and epithelium. In the regime of small strains, the coupling mechanisms include degradation of FGF due to local solid compression, advection of chemicals through solid velocity, and an active stress driven by mesenchymal cells. We have derived a set of conditions for Turing patterns to emerge, focusing on the two-way coupling and the main coupling variables. We have proposed an extended model that accounts for growth and hyperelastic deformations of the solid matrix, and the coupling mechanisms in this case consist in intrinsic growth functions depending on time and on the concentration of mesenchymal cells, as well as modification of the diffusion coefficients in the chemotaxis model due to changes from spatial to reference coordinates.

The linear stability analysis for the nonlinearly coupled system from Section~\ref{sec:growth} could be carried out extending our results from Section~\ref{sec:stability} following the ideas from \cite{benamar05}. As the main application is on skin patterning for feather development, we could also incorporate growth models specifically targeted for shells or thin plates as in \cite{dervaux09}, and  concentrate on bi-layered structures extending the formulations in \cite{deoliveira19}.

\section*{Acknowledgements}
This research was funded by  the Monash Mathematics Research Fund S05802-3951284; and by the Australian Research Council through the Discovery Project grant DP220103160.
\section*{Appendix}

\appendix

\subsection*{Proof of Proposition~\ref{prop:ImpDerivative}.}
The spatial derivatives in the scalar case follow trivially from the definition of $f$. Both scalar and vector time derivatives also follow directly from the definition of $f,\bv$. For the spatial derivatives of $\bv$, using the definition of the differential operators and denoting $\{\ev_n\}_{1,\ldots,d}$, the standard basis vectors in $\mathbb{R}^d$, we can readily write
\begin{linenomath*}	\begin{align*}
		\vdiv \bv &= \sum_n \partial_{x_n} v_{n,0} \e^{\ii \bk\cdot\bx + \phi t} = \sum_n  \ii k_n v_{n,0} \e^{\ii \bk\cdot\bx + \phi t} = \sum_n \ii k_n v_n = \ii \bv \cdot \bk, \\
		\bnabla \bv &= \sum_m \partial_{x_m} \bv \otimes \ev_m = \sum_{n,m} \partial_{x_m} v_{n,0} \e^{\ii \bk\cdot\bx + \phi t} \ev_n \otimes \ev_m = \sum_{n,m} \ii k_m v_{n,0} \e^{\ii \bk\cdot\bx + \phi t} \ev_n \otimes \ev_m \\
		&= \sum_{n,m} \ii k_m v_n \ev_n \otimes \ev_m = \sum_m \ii k_m \bv \otimes \ev_m = \ii \bv \otimes \bk.
	\end{align*}\end{linenomath*}
On the other hand, using \eqref{eq:DerivativeFV} we can then straightforwardly derive  \eqref{eq:DerivativeEpsilonTheta} as follows
\begin{linenomath*}\begin{align*}
		\bdiv \beps(\bv) &= \bdiv \left( \sum_{n,m} \ii \frac{v_nk_m + v_m k_n}{2} \ev_n \otimes \ev_m \right)
		= \sum_{n,m} \ii^2 \frac{v_nk_nk_m + v_m k_n^2}{2} \ev_m \\
		& = -\frac{1}{2} \sum_m (\bv \cdot \bk)k_m \ev_m + k^2 v_m \ev_m = -\frac{(\bv\cdot\bk)\bk + k^2\bv}{2}, \\
	\partial_t \beps(\bv) &= \left( \frac{1}{2}(\bnabla \partial_t\bv + \bnabla \partial_t\bv^\intercal) \right)= \left( \frac{\phi}{2}(\bnabla \bv + \bnabla \bv^\intercal) \right) = \phi \beps(\bv),  \\
		\bdiv (\theta \bI) &= \nabla\theta = \nabla(\ii\bv\cdot\bk) = \ii \sum_n \partial_{x_n} (\bv \cdot \bk) \ev_n = \ii \sum_{n,m} k_m \partial_{x_n}v_m \ev_n = \ii^2 \sum_{n,m} k_m k_n v_m \ev_n  \\
		&= -\sum_n k_n (\bv\cdot\bk)\ev_n = - (\bv\cdot\bk)\bk, \\
		\partial_t \theta  &= \vdiv \partial_t \bv = \phi \vdiv \bv = \phi \theta.
	\end{align*}\end{linenomath*}

\end{document}